# LOCAL LIMIT OF LABELED TREES AND EXPECTED VOLUME GROWTH IN A RANDOM QUADRANGULATION

By Philippe Chassaing and Bergfinnur Durhuus[1]

*Université H. Poincaré–Nancy and University of Copenhagen*

Exploiting a bijective correspondence between planar quadrangulations and well-labeled trees, we define an ensemble of infinite surfaces as a limit of uniformly distributed ensembles of quadrangulations of fixed finite volume. The limit random surface can be described in terms of a birth and death process and a sequence of multitype Galton–Watson trees. As a consequence, we find that the expected volume of the ball of radius $r$ around a marked point in the limit random surface is $\Theta(r^4)$.

**1. Introduction.** Beginning with the seminal article [29] by Tutte, the combinatorics of planar maps have been a subject of continuing development, gaining further impetus in recent years after the realization of its importance in quantum field theory [11] and in string theory and two-dimensional quantum gravity [3, 14, 21]. In the latter case planar maps play the role of two-dimensional discretized (Euclidean) space–time manifolds, whose topology equals that of a sphere with a number of holes; see, for example, [5] for an overview. We shall frequently use the term "planar surface" instead of "planar map."

Of primary interest up to now have been properties that depend on the volume, that is, the number of vertices, the number and length of the boundary components or local properties such as the distribution of vertex degrees, whereas quantities that depend on the internal (graph) metric structure have received less attention, despite their obvious relevance, for example, in quantum gravity [4]. Given a probability measure on a space of planar maps (in physics, also called an *ensemble* of planar maps), a relevant quantity to

Received December 2003; revised October 2004.
[1]Supported in part by MaPhySto—A Network in Mathematical Physics and Stochastics, funded by the Danish National Research Foundation.

consider in this connection is the exponent $\alpha$ that characterizes their expected volume growth. In the physics literature (see, e.g., [4]), this exponent is sometimes called the *Hausdorff dimension* of the ensemble of surfaces. It is defined more precisely as follows. Denoting by $B_r(\mathcal{M})$ a ball of radius $r$ around a marked point in the surface $\mathcal{M}$ and denoting by $|B_r(\mathcal{M})|$ its volume, that is, the number of vertices in $B_r(\mathcal{M})$, then $\alpha$ is determined by

$$\mathbb{E}[|B_r|] = \Theta(r^\alpha), \tag{1}$$

assuming such a relationship exists. Here the expectation value is understood with respect to the probability measure in question, and we use the standard notation $\Theta(\phi(r))$ for a function bounded from above and below by positive constant multiples of $\phi(r)$ for $r$ large enough. In this paper we use the notation $\mathbb{E}[\cdots]$ rather than $\langle\cdots\rangle$ for the expectation of a random variable. It is implicit in the definition that the support of the measure consists of surfaces of infinite extent. We shall define such a measure as a limit of uniformly distributed surfaces of fixed finite volume.

It is also possible, and more common, to define $\alpha$ directly in terms of an ensemble of finite surfaces. For instance, one can introduce the so-called two-point correlation function $G(r)$, which is defined as a certain integral over finite planar surfaces with two marked points a fixed (graph) distance $r$ apart, whose exponential decay rate as a function of $r$ can be shown to determine $\alpha$; see [5]. This method was applied to triangulated planar surfaces (or pure two-dimensional quantum gravity) in [6], where it was argued that $\alpha = 4$ (see also [10]).

A different method was recently applied in [12] to quadrangulated planar maps, making it possible to exploit a clever bijective correspondence between such maps and so-called *well-labeled trees*. That work established the existence of the limit in law of the random variables $(n^{-1}|B_{xn^{1/4}}|)_{x\geq 0}$ and $n^{-1/4}r_n$ as $n \to \infty$, where $r_n$ (resp. $|B_r|$) is the radius (resp. the volume of a ball) of a random quadrangulation with $n$ faces, identifying the limits as functionals of the so-called one-dimensional integrated super-Brownian excursion (or ISE) [2, 28]. In view of the information carried by these results, it seems reasonable to make the identification $\alpha = 4$.

Yet another approach was recently proposed in [7, 8]. In [8] a *uniform probability measure* is first constructed on *infinite planar triangulations* as a limit of uniform measures on finite planar triangulations. In [7] it is then proven that the ball of radius $r$ around a marked triangle in a sample triangulation contains on the order of $r^4$ triangles up to logarithmic corrections, which can be seen as yet another manifestation of the exponent $\alpha$ being equal to 4.

In this article we adapt the technique of [8] to the case of well-labeled trees and thereby construct by simple combinatorial arguments a uniform probability measure $\mu$ on infinite well-labeled trees. We show how to identify well-labeled trees in the support of this measure with infinite quadrangulated planar surfaces through a mapping that shares the basic properties of Schaeffer's bijection [12, 27] between finite trees and finite quadrangulations. In particular, there is a one-to-one correspondence between the vertices in a tree and the vertices in the corresponding surface, except for a certain marked vertex in the surface, and the label $r \in \mathbb{N}$ of a vertex in a tree equals the (graph) distance between the corresponding vertex in the surface and the marked vertex.

Via this identification, viewing $\mu$ as a measure on quadrangulated planar surfaces, we prove the relationship (1) with $\alpha = 4$. Thus the exponent 4 does not appear only for the volume of large balls in large but finite random quadrangulations (cf. [12]). The result of the present paper, for balls with radius $r$ in an ensemble of infinite quadrangulations, is closer in spirit to the result of [7] for triangulations. Together these results corroborate the still unproven claim of universality of the exponent $\alpha = 4$ for planar random surfaces.

The article is organized as follows. In Section 2 we define the space of rooted (labeled) trees. It is endowed with a topology, and its topological and combinatorial properties are studied. In Section 3 we construct a so-called uniform measure on the set of infinite rooted trees, whose vertices are labeled by positive integers, such that the labels of neighboring vertices deviate by at most 1 and such that the root has a fixed label. In case the root has label 1, such trees are called *well-labeled* in [12]. In Section 4 we show that almost surely the trees have exactly one infinite branch, allowing a definition of the spine of a sample well-labeled tree as the unique infinite non-self-intersecting path starting at the root. A sample tree can be obtained by attaching (finite) branches along the spine, independently distributed according to measures $\hat{\rho}^{(k)}$, $k$ denoting the label of the root of a branch. In particular, we show in Section 4.2 that the labels along the spine are described by a certain birth and death process, and that the labels in the branches are described in terms of a multitype Galton–Watson process. In Section 5 these two processes are investigated in more detail. In particular, the birth and death process is shown to be transient and we determine, as a main result, the asymptotic behavior of the average number $\mathbb{E}[N_r]$ of vertices with a fixed label $r$ to be

$$(2) \qquad \mathbb{E}[N_r] = \Theta(r^3)$$

for $r$ large. In Section 6 we show how to extend the mapping of well-labeled trees onto planar quadrangulations to infinite ones. Combined with (2) this yields (1) and the value 4 for the exponent $\alpha$ of planar quadrangulated surfaces.

**2. Labeled trees.**

2.1. *Definitions and notation.* By $\overline{\mathcal{T}}$ we shall denote the set of rooted planar trees, where *rooted* means that one oriented edge $(i_0, i_1)$ is distinguished, called the *root*, and $i_0$ and $i_1$ are called the first and second *root vertex*, respectively. Here, trees are allowed to be infinite, but vertices are of finite degree. The adjective planar means that trees are assumed to be embedded into the plane $\mathbb{R}^2$ such that no two edges intersect except at common vertices, and we identify trees that can be mapped onto each other by an orientation-preserving homeomorphism of the plane that maps root onto root. In addition, certain regularity requirements on the embeddings are needed, the discussion of which we postpone until Section 6. A more precise combinatorial definition is as follows. Once a fixed orientation of $\mathbb{R}^2$ is chosen, the set of vertices at distance $r$ from the first root vertex in a given rooted planar tree $\tau$ has a natural ordering. This can be obtained, for example, by choosing a right-handed coordinate system for $\mathbb{R}^2$ and mapping the tree into $\mathbb{R}^2$ such that the vertices at distance $r$ from the root $i_0$ are mapped into the the vertical line through $(r, 0)$ and then ordering according to their second coordinate in such a way that, for $r = 1$, the second root $i_1$ is smallest. We call this ordered set $\Phi_r = (i_{r1}, \ldots, i_{rn_r})$. The edges in the tree are specified by mappings $\phi_r : \Phi_r \to \Phi_{r-1}$, $r \geq 1$, that preserve the ordering, that is, the edges in $\tau$ are given by $(i_{rk}, i_{r-1\phi_r(k)})$, $1 \leq k \leq n_r$. It is clear that any (finite or infinite) sequence $(\Phi_0, \Phi_1, \Phi_2, \ldots)$ of finite pairwise disjoint ordered sets, where $\Phi_0 = i_0$ is a one-point set, together with orientation-preserving maps $(\phi_1, \phi_2, \ldots)$ as above, uniquely specifies a rooted planar tree $\tau$, in which case we write $\tau = (\Phi_r, \phi_r)_{r \in \mathbb{N}}$. We then have $(\Phi_r, \phi_r)_{r \in \mathbb{N}} = (\Phi'_r, \phi'_r)_{r \in \mathbb{N}}$ if and only if there exist order-preserving bijective maps $\psi_r : \Phi_r \to \Phi'_r$ such that $\phi'_r = \psi_{r-1} \circ \phi_r \circ \psi_r^{-1}$ for all $r$.

We have

$$\overline{\mathcal{T}} = \left( \bigcup_{N=1}^{\infty} \mathcal{T}_N \right) \cup \mathcal{T}_\infty,$$

where $\mathcal{T}_N$ consists of trees with maximal vertex distance from the first root equal to $N$, that is, $\Phi_r = \varnothing$ for $r > N$ but $\Phi_N \neq \varnothing$, and $\mathcal{T}_\infty$ consists of infinite trees. We say that $\tau \in \mathcal{T}_N$ has *height* (or radius) $\rho(\tau) = N$. The set $\bigcup_{N=1}^{\infty} \mathcal{T}_N$ of finite trees will be denoted by $\mathcal{T}$ and the *size* (or volume) $|\tau|$ of a finite tree is defined to be the number of edges in $\tau$. For $\tau \in \mathcal{T}_\infty$ we set $\rho(\tau) = |\tau| = \infty$.

By a labeled tree we mean a pair $(\tau, \ell)$, where $\ell : i \to \ell_i$ is a mapping from the vertices of $\tau$ into the integers $\mathbb{Z}$, such that

$$|\ell_i - \ell_j| \leq 1 \quad \text{if } (i, j) \text{ is an edge in } \tau.$$

If, furthermore,

$$\ell_{i_0} = k \quad \text{and} \quad \ell_i \geq 1 \quad \text{for all vertices } i \text{ in } \tau,$$

we say that $(\tau, \ell)$ is a *k-labeled* tree. A 1-labeled tree is also called a *well-labeled tree* [12, 13]. We call the set of $k$-labeled trees $\overline{\mathcal{W}}^{(k)}$ and, for $k = 1$, we set $\overline{\mathcal{W}}^{(1)} = \overline{\mathcal{W}}$. The corresponding sets of finite labeled trees are denoted similarly without overbars. Obviously, $\overline{\mathcal{W}}^{(k)}$ and $\overline{\mathcal{W}}$ inherit from $\overline{\mathcal{T}}$ the natural decompositions

$$\overline{\mathcal{W}}^{(k)} = \mathcal{W}^{(k)} \cup \mathcal{W}^{(k)}_\infty = \left(\bigcup_{N=1}^\infty \mathcal{W}^{(k)}_N\right) \cup \mathcal{W}^{(k)}_\infty,$$

$$\overline{\mathcal{W}} = \mathcal{W} \cup \mathcal{W}_\infty = \left(\bigcup_{N=1}^\infty \mathcal{W}_N\right) \cup \mathcal{W}_\infty,$$

into finite and infinite $k$-labeled trees. If $\omega = (\tau, \ell)$ is a labeled tree, we set $|\omega| = |\tau|$ and $\rho(\omega) = \rho(\tau)$. Moreover, if $\tau = (\Phi_r, \phi_r)_{r \in \mathbb{N}}$ and $\Phi_r = (i_{r1}, \ldots, i_{rn_r})$, we set

$$\Xi_r = ((i_{r1}, \ell_{i_{r1}}), \ldots, (i_{rn_r}, \ell_{i_{rn_r}})),$$

in which case we have

(3) $\qquad |\ell_{rk} - \ell_{r-1\phi_i(k)}| \leq 1 \qquad$ for all $r \geq 1, 1 \leq k \leq n_r$.

Clearly, any sequence $(\Xi_r, \phi_r)$, $r \in \mathbb{N}$, obtained from a tree $\tau = (\Phi_r, \phi_r)_{r \in \mathbb{N}}$ by adding labels as above to the vertices of each $\Phi_r$ fulfilling (3) defines a unique labeled tree $\omega$, in which case we write $\omega = (\Xi_r, \phi_r)_{r \in \mathbb{N}}$. Furthermore, we have $(\Xi_r, \phi_r)_{r \in \mathbb{N}} = (\Xi'_r, \phi'_r)_{r \in \mathbb{N}}$ if and only if there exist maps $\psi_r$ that identify the underlying unlabeled trees and such that the labels of a vertex $i$ and its image $\psi_r(i)$ are identical.

2.2. *Topology on labeled trees.* For $r \in \mathbb{N}_0$ and $\omega \in \mathcal{W}^{(k)}$ with distance classes $\Xi_s(\omega)$, $s \in \mathbb{N}_0$, we define the ball $B_r(\omega)$ of radius $r$ in $\omega$ to be the labeled subtree of $\omega$ generated by $\Xi_0(\tau), \ldots, \Xi_r(\tau)$ if $r < \rho(\omega)$ and equal to $\omega$ otherwise. In other words,

$$B_r(\omega) = ((\Xi_1, \phi_1), \ldots, (\Xi_r, \phi_r)) \qquad \text{if } \omega = ((\Xi_1, \phi_1), (\Xi_2, \phi_2), \ldots).$$

Next we define, for $\omega, \omega' \in \overline{\mathcal{W}}^{(k)}$ and $k$ fixed,

$$d(\omega, \omega') = \inf\left\{\frac{1}{r+1}\bigg| B_r(\omega) = B_r(\omega'), r \in \mathbb{N}_0\right\}.$$

It is trivially verified that $d$ defines a metric on $\overline{\mathcal{W}}^{(k)}$. The corresponding open balls in $\overline{\mathcal{W}}^{(k)}$ are given by

$$\mathcal{B}_s(\omega_0) = \{\omega \in \overline{\mathcal{W}}^{(k)} | d(\omega, \omega_0) < s\} \qquad \text{for } s > 0.$$

REMARK 2.1. The following facts are easy to verify:

- The set $\mathcal{W}^{(k)}$ of finite $k$-labeled trees is a countable dense subset of $\overline{\mathcal{W}}^{(k)}$ and its boundary $\partial \mathcal{W}^{(k)}$ in $\overline{\mathcal{W}}^{(k)}$ equals $\mathcal{W}_\infty^{(k)}$.
- For $s > 0$ and $\omega \in \overline{\mathcal{W}}^{(k)}$, the ball $\mathcal{B}_s(\omega)$ is both open and closed, and

$$\{\omega \in \mathcal{B}_s(\omega_0)\} \Leftrightarrow \{\mathcal{B}_s(\omega) = \mathcal{B}_s(\omega_0)\}.$$

As a consequence, either two balls are disjoint or one is contained in the other.

- The set $\overline{\mathcal{W}}^{(k)}$ is not compact: Let $\omega_n$ be the (unique) $k$-labeled tree of height 1 with $n+1$ vertices and all labels equal to $k$. Then $d(\omega_n, \omega_m) = 1$ for $n \neq m$ and hence $\omega_n, n \in \mathbb{N}$, has no convergent subsequence.

As a substitute for compactness, we shall make use of the following result:

PROPOSITION 2.2. *Let $K_r, r \in \mathbb{N}$, be a sequence of positive numbers. Then the subset*

$$C = \bigcap_{r=1}^\infty \{\omega \in \overline{\mathcal{W}}^{(k)} | |B_r(\omega)| \leq K_r\}$$

*of $\overline{\mathcal{W}}^{(k)}$ is compact.*

PROOF. Let $\omega_n, n \in \mathbb{N}$, be any sequence in $C$. For each $r \in \mathbb{N}$, the set

$$\{\omega \in \mathcal{W}_r^{(k)} | |\omega| \leq K_r\}$$

is finite. Hence there exists a subsequence $\omega_{n_i}, i \in \mathbb{N}$, such that $B_r(\omega_{n_i})$ is constant as a function of $i$. By applying a diagonal argument, we may choose this subsequence such that $B_i(\omega_{n_j}) = B_i(\omega_{n_i})$ for all $i \leq j$. It follows that this subsequence determines a unique tree $\omega \in C$ such that $B_i(\omega) = B_i(\omega_{n_i})$ for all $i \in \mathbb{N}$. In particular, $\omega_{n_i} \to \omega$ as $i \to \infty$, which completes the proof. □

2.3. *Combinatorics of finite labeled trees.* In Section 3 we shall consider the sequence $\mu_N$, $N \in \mathbb{N}$, of measures on $\overline{\mathcal{W}}$, where $\mu_N$ is defined as the uniform probability measure concentrated on

$$\mathcal{W}'_N = \{\omega \in \mathcal{W} | |\omega| = N\}, \qquad N \in \mathbb{N}_0,$$

that is,

$$\mu_N(\omega) = D_N^{-1} \qquad \text{for } \omega \in \mathcal{W}'_N, \mu_N(\mathcal{W} \setminus \mathcal{W}'_N) = 0,$$

where $D_N = \sharp \mathcal{W}'_N$ is the number of well-labeled trees of size $N$.

To establish weak convergence of $\mu_N$, $N \in \mathbb{N}$, we need some basic facts about the sequence $D_N$, $N \in \mathbb{N}$, and, more generally, about the sequence $D_N^{(k)}$, $N \in \mathbb{N}$, where $D_N^{(k)}$ is the number of $k$-labeled trees of size $N$.

As shown in [12, 13], $D_N$ equals the number of quadrangulated planar maps with $N$ faces; see also Section 6 below. The corresponding generating function $W(z)$ has been computed in [29] and is given by

$$
\begin{align}
W(z) &= \sum_{N=0}^{\infty} D_N z^N \\
&= 1 + 2z + 9z^2 + \cdots \\
&= \frac{18z - 1 + (1 - 12z)^{3/2}}{54z^2} \qquad \text{for } |z| \leq \frac{1}{12},
\end{align}
\tag{4}
$$

yielding

$$
D_N = 2 \cdot 3^N \frac{(2N)!}{N!(N+2)!}.
\tag{5}
$$

Note that we have included in $W(z)$ the contribution $D_0 = 1$ from the tree with only one vertex. As consequences of (5) we have

$$
D_N \simeq \frac{2}{\sqrt{\pi}} N^{-5/2} 12^N
\tag{6}
$$

and, for each $k \in \mathbb{N}$,

$$
D_N^{(k)} = \Theta(N^{-5/2} 12^N),
\tag{7}
$$

relationship (7) being a consequence of

$$
D_{N-k+1} \leq D_N^{(k)} \leq D_{N+k-1}.
$$

These last inequalities can be obtained by grafting a suitable tree on a branch with length $k$ and with labels strictly increasing or decreasing from 1 to $k$.

We need a more precise estimate for $D_N^{(k)}$. For this purpose, we define

$$
W^{(k)}(z) = \sum_{N=0}^{\infty} D_N^{(k)} z^N,
$$

$$
Z^{(k)}(z) = \sum_{N=1}^{\infty} E_N^{(k)} z^N,
$$

where $E_N^{(k)}$ is the number of $k$-labeled trees with $N$ edges and first root $i_0$ of degree 1 and we have set $D_0^{(k)} = 1$. Note that estimates like (7) hold true

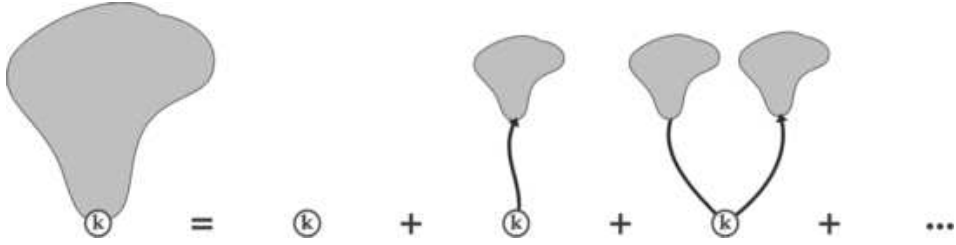

Fig. 1.  *Decomposing a tree:* $W^{(k)} = \frac{1}{1-Z^{(k)}}$.

also for $E_N^{(k)}$, by a similar proof. By these estimates, $W^{(k)}(z)$ and $Z^{(k)}(z)$ are analytic in the open disc

$$\Delta_0 = \{z \in \mathbb{C} | |z| < \tfrac{1}{12}\}$$

and extend to continuous functions on its closure $\bar\Delta_0$. By decomposing a tree into trees with first root $i_0$ of degree 1 as in Figure 1, we obtain

LEMMA 2.3. *In $\bar\Delta_0$, $W^{(k)}$ and $1 - Z^{(k)}$ have no roots, and*

$$(8) \qquad W^{(k)} = \frac{1}{1 - Z^{(k)}}.$$

In the following, we shall need the values

$$z_k = Z^{(k)}(\tfrac{1}{12}) \quad \text{and} \quad w_k = W^{(k)}(\tfrac{1}{12}).$$

PROPOSITION 2.4. *For $k \in \mathbb{N}$ we have*

$$(9) \qquad z_k = \frac{1}{2} - \frac{1}{k(k+3)}$$

*and*

$$(10) \qquad w_k = 2\frac{k(k+3)}{(k+1)(k+2)}.$$

PROOF. By (8) it suffices to prove (9). For this purpose we note the relationships

$$(11) \quad Z^{(1)}(z) = \frac{z}{1 - Z^{(1)}(z)} + \frac{z}{1 - Z^{(2)}(z)},$$

$$(12) \quad Z^{(k)}(z) = \frac{z}{1 - Z^{(k-1)}(z)} + \frac{z}{1 - Z^{(k)}(z)} + \frac{z}{1 - Z^{(k+1)}(z)}, \qquad k \geq 2,$$

which are obtained by decomposing the sum over trees that define $Z^{(k)}$ according to the degree and the label of their second root $i_1$.

Clearly, these relationships determine $Z^{(k)}$ in terms of $Z^{(1)}$. Inserting $z = \frac{1}{12}$ and the value $z_1 = \frac{1}{4}$ obtained from (4), one finds that (9) solves (11) and (12). □

A precise estimate for $D_N^{(k)}$ is established in the subsequent lemma. Set

$$d_{k,N} = \frac{D_N^{(k)}}{D_N}. \tag{13}$$

LEMMA 2.5. *For each $k \in \mathbb{N}$, the sequence $(d_{k,N})_{N \geq 1}$ converges to a limit*

$$d_k = \frac{3}{280} \frac{k(k+3)}{(k+1)(k+2)} (5k^4 + 30k^3 + 59k^2 + 42k + 4). \tag{14}$$

*This limit fulfills*

$$d_1 = 1, \tag{15}$$

$$d_1 + d_2 = 12 d_1 (w_1)^{-2}, \tag{16}$$

$$d_{k-1} + d_k + d_{k+1} = 12 d_k (w_k)^{-2}, \qquad k \geq 2. \tag{17}$$

From (14) it follows that $d_k$ is positive, so that we have the following corollary:

COROLLARY 2.6. *We have*

$$D_N^{(k)} \simeq \frac{2 d_k}{\sqrt{\pi}} N^{-5/2} 12^N.$$

PROOF OF LEMMA 2.5. Due to expression (4), $W = W^{(1)}$ extends to an analytic function on $\mathbb{C}$, except for a cut along the positive real axis starting at $z = \frac{1}{12}$. In particular, $W^{(1)}$ is analytic on the set

$$\Delta_r = \{z \in \mathbb{C} | |z| < \tfrac{1}{12} + r\} \setminus \{\tfrac{1}{12} + s e^{i\theta} | s \geq 0, |\theta| \leq \epsilon\}$$

for any $r, \epsilon > 0$ and is continuous on $\bar{\Delta}_r$. We henceforth fix $0 < \epsilon < \frac{\pi}{2}$. By Lemma 2.3 and uniform continuity on $\bar{\Delta}_r$, we can choose $r(1) > 0$ such that $W^{(1)}(z) \neq 0$ on $\bar{\Delta}_{r(1)}$. Rewriting relationships (8), (11) and (12) in the form

$$W^{(k+1)}(z) = \frac{1}{z}\left(1 - \frac{1}{W^{(k)}(z)}\right) - W^{(k-1)}(z) - W^{(k)}(z), \qquad k \geq 1, \tag{18}$$

where we use the convention $W^{(0)} = 0$, it follows that $W^{(2)}$ extends to a continuous function on $\bar{\Delta}_{r(1)}$ that is analytic on $\Delta_{r(1)}$. Again $r(2) > 0$ can be chosen small enough such that $W^{(2)}(z) \neq 0$ on $\bar{\Delta}_{r(2)}$. By induction we conclude that each $W^{(k)}$ is analytic on a domain $\Delta_{r(k)}$, $r(k)$ positive and extends continuously to its closure.

In particular, $\Delta_{r(k)}$ is a $\Delta$ domain for $W^{(k)}$ in the sense of [17] or [18], Chapter 6.3. Furthermore, note that an expansion of the form

$$W^{(k)}(z) = w_k + a_k(1 - 12z) + b_k(1 - 12z)^{3/2} + o((1 - 12z)^{3/2}),$$

which holds for $k = 1$ due to (4), also holds for general $k$ as a consequence of (18). Identification of the coefficients on the two sides of (18) yields the relationship

$$(19) \qquad b_{k+1} = \frac{12b_k}{w_k^2} - b_{k-1} - b_k, \qquad k \geq 1,$$

with the convention $b_0 = 0$. Thus, transfer theorems around $1/12$ (cf. [18], Chapter 6.3, Theorems 5.4 and 5.5) imply that

$$\lim_N D_N^{(k)} N^{5/2} 12^{-N} = \frac{3b_k}{4\sqrt{\pi}}.$$

It follows that $d_k$ exists and

$$d_k = \frac{b_k}{b_1} = \frac{3b_k}{8}.$$

Relationships (16) and (17) are then immediate consequences of (19). Since $D_N = D_N^{(1)}$, the limit $d_1$ trivially exists and equals 1. It turns out (tediously checking identities between polynomials of order 8) that the expression given in (14) satisfies (15), (16) and (17). Clearly, this last set of equations determines $d_k$ uniquely. $\square$

REMARK 2.7. In [10] (cf. [10], (4.19)), it is mentioned that (14) can be derived from a closed form expression of $W^{(k)}$ by a saddle-point approximation. Above we proposed an alternative proof.

**3. Uniform measure on infinite labeled trees.** We are now ready to prove one of the main results of this paper.

THEOREM 3.1. *The sequence $\mu_N$, $N \in \mathbb{N}$, converges weakly to a Borel probability measure, $\mu$, concentrated on $\mathcal{W}_\infty$. We call $\mu$ the uniform probability measure on $\mathcal{W}_\infty$.*

PROOF. By Remark 2.1 the denumerable family of balls

$$\mathcal{V} = \{\mathcal{B}_{1/r}(\omega) | r \in \mathbb{N}, |\omega| < +\infty\}$$

consists of open and closed sets. In addition:

(i) Any finite nonempty intersection of sets in $\mathcal{V}$ belongs to $\mathcal{V}$.
(ii) Any open set in $\overline{\mathcal{W}}$ can be written as a union of sets in $\mathcal{V}$.

By Theorems 2.1, 2.2 and 6.1 in [9] it will suffice to prove that the sequence $\mu_N$, $N \in \mathbb{N}$, is tight and that $\mu_N(A)$ converges as $N \to \infty$ for all $A \in \mathcal{V}$.

We first prove tightness by showing that for any given $\varepsilon > 0$ and $r \in \mathbb{N}$ there exists $K_r > 0$ such that

$$\mu_N(\{\omega \in \mathcal{W} | |B_r(\omega)| > K_r\}) < \varepsilon \tag{20}$$

for all $N$. Replacing $\varepsilon$ by $\varepsilon/2^r$ in (20) and choosing $K_r$ correspondingly, Proposition 2.2 gives the desired compact set $C$ fulfilling $\mu_N(C) > 1 - \varepsilon$ for all $N$.

We proceed to show (20) by induction on $r$. If $r = 1$, then $|B_r(\omega)|$ equals the degree of the root vertex $i_0(\omega)$. By the argument leading to (8) there is a one-to-one correspondence between trees $\omega$ in $\mathcal{W}_N$ with $i_0(\omega)$ of degree $K$ and $K$-tuples $(\omega_1, \ldots, \omega_K)$ of trees, such that the first root vertex of each $\omega_a$ has degree 1 and

$$|\omega_1| + \cdots + |\omega_K| = N.$$

This gives

$$\mu_N(\{\omega \in \mathcal{W} | |B_1(\omega)| = K\}) = D_N^{-1} \sum_{N_1 + \cdots + N_K = N} \prod_{a=1}^{K} E_{N_a}^{(1)}.$$

In this sum, $N_a \geq N/K$ for at least one value of $a = 1, \ldots, K$. Combining this with (7) and

$$E_N^{(k)} \leq D_N^{(k)},$$

we obtain

$$\mu_N(\{\omega \in \mathcal{W} | |B_1(\omega)| = K\}) \leq K \sum_{\substack{N_1 + \cdots + N_K = N \\ N_1 \geq N/K}} c' \cdot K^{5/2} \prod_{a=2}^{K} E_{N_a}^{(1)} \cdot 12^{-N_a}$$

$$\leq c' \cdot K^{7/2} z_1^{K-1} = c' K^{7/2} 4^{1-K},$$

where $c' > 0$ is a constant independent of $N$. This proves (20) for $r = 1$ and for sufficiently large $K_1$.

Now assume (20) holds for a given $r \geq 1$. For any $K > 0$, we then have

$$\mu_N(\{\omega \in \mathcal{W} | |B_{r+1}(\omega)| > K\})$$
$$\leq \varepsilon + \mu_N(\{\omega \in \mathcal{W} | |B_{r+1}(\omega)| > K, |B_r(\omega)| \leq K_r\}).$$

Since there are only finitely many different balls $B_r(\omega)$ with $|B_r(\omega)| \leq K_r$, it suffices to show that

$$\mu_N(\{\omega \in \mathcal{W} | |B_{r+1}(\omega)| > K, B_r(\omega) = \hat{\omega}\}) \to 0 \tag{21}$$

as $K \to \infty$, uniformly in $N$ for any fixed $\hat{\omega} \in \mathcal{W}_r$. This is obtained in a similar fashion as for $r = 1$:

Set $\Xi_r(\hat{\omega}) = ((j_1, k_1), \ldots, (j_R, k_R))$ and $K' = K - |\hat{\omega}|$. For $\omega \in \mathcal{W}$ with $B_r(\omega) = \hat{\omega}$ and $|B_{r+1}(\omega)| = K$, let $(\omega_1, \ldots, \omega_{K'})$ be the ordered set of branches with roots of degree 1 attached to the vertices of $\Xi_r(\hat{\omega})$, such that the first $L_1 \geq 0$ are attached to $j_1$, the next $L_2 \geq 0$ are attached to $j_2, \ldots$ and the last $L_R$ are attached to $j_R$. Then

$$L_1 + \cdots + L_R = K' \tag{22}$$

and the labels $\ell_1, \ldots, \ell_{K'}$ of the first root vertices of $\omega_1, \ldots, \omega_{K'}$ are determined by $k_1, \ldots, k_R$. Hence, the subset of trees $\omega$ with fixed values of $L_1, \ldots, L_R \geq 0$ fulfilling (22) has $\mu_N$ measure given by

$$D_N^{-1} \sum_{N_1+\cdots+N_{K'}=N-|\hat{\omega}|} \prod_{s=1}^{K'} E_{N_s}^{(\ell_s)} \leq D_N^{-1} \sum_{t=1}^{K'} \sum_{\substack{N_1+\cdots+N_{K'}=N-|\hat{\omega}| \\ N_t \geq (N-|\hat{\omega}|)/K'}} \prod_{s=1}^{K'} E_{N_s}^{(\ell_s)}$$

$$\leq c'' \sum_{t=1}^{K'} K'^{5/2} \prod_{\substack{s=1 \\ s \neq t}}^{K'} z_{\ell_s}$$

$$\leq c'' K^{7/2} 2^{|\hat{\omega}|+1-K},$$

where the first two inequalities follow by the same arguments as for $r = 1$ and where $c''$ is a constant that depends only on $\hat{\omega}$. Since there are

$$\binom{K' + R - 1}{R - 1} \leq \frac{K^{R-1}}{(R-1)!}$$

different ways to choose $L_1, \ldots, L_R \geq 0$ that fulfill (22), we have

$$\mu_N(\{\omega \in \mathcal{W} | |B_{r+1}(\omega)| = K, B_r(\omega) = \hat{\omega}\}) \leq c''' K^{R+5/2} 2^{-K}$$

for some positive constant $c'''$ that depends only on $\hat{\omega}$. Clearly, this proves (21) for $r + 1$ if $K_{r+1}$ is chosen large enough. This completes the proof of tightness of the sequence $\mu_N$, $N \in \mathbb{N}$.

It remains to establish convergence of $\mu_N(\mathcal{B}_{1/r}(\hat{\omega}))$ as $N \to \infty$ for all $r \in \mathbb{N}$ and finite $\hat{\omega} \in \mathcal{W}$. When $\rho(\hat{\omega}) \leq r - 1$, we have $\mathcal{B}_{1/r}(\hat{\omega}) = \{\hat{\omega}\}$ and

$$\lim_N \mu_N(\mathcal{B}_{1/r}(\hat{\omega})) = 0.$$

Since $\mathcal{B}_{1/r}(\hat{\omega}) = \mathcal{B}_{1/r}(B_r(\hat{\omega}))$, we assume from now on that $\hat{\omega} \in \mathcal{W}_r$ and we set $\hat{N} = |\hat{\omega}|$. Then

$$\mathcal{B}_{1/r}(\hat{\omega}) = \{\omega \in \mathcal{W} | B_r(\omega) = \hat{\omega}\}.$$

With the notation $\Delta_r(\hat{\omega}) = (j_1, \ldots, j_R)$, any $\omega \in \mathcal{B}_{1/r}(\hat{\omega})$ is obtained by grafting a sequence $(\omega_1, \ldots, \omega_R)$ of $R$ trees in $\mathcal{W}$ on $\hat{\omega}$ such that the root vertex $i_0(\omega_s)$ is identified with $j_s$ and has the same label $k_s$ as that of $j_s$ in $\hat{\omega}$. This gives

$$(23) \quad \mu_N(\mathcal{B}_{1/r}(\hat{\omega})) = D_N^{-1} \sum_{N_1 + \cdots + N_R = N - \hat{N}} \prod_{s=1}^{R} D_{N_s}^{(k_s)},$$

where $N_s = |\omega_s|$. For a given $t = 1, \ldots, R$, let $N_s, s \neq t$, be fixed, while $N_t$ is determined as a function of $N$ by $N_1 + \cdots + N_R = N - \hat{N}$. For the corresponding term in the sum, we obtain from Corollary 2.6, as $N \to \infty$,

$$D_N^{-1} \prod_{s=1}^{R} D_{N_s}^{(k_s)} = 12^{-\hat{N}} \frac{D_{N_t}^{(k_t)}}{D_{N_t}} \frac{12^{-N_t} D_{N_t}}{12^{-N} D_N} \prod_{s \neq t} D_{N_s}^{(k_s)} 12^{-N_s}$$

$$\to 12^{-\hat{N}} d_{k_t} \prod_{s \neq t} D_{N_s}^{(k_s)} 12^{-N_s}.$$

From this we conclude, for any fixed $A > 0$, that

$$(24) \quad \lim_N \mu_N(\{\omega \in \mathcal{B}_{1/r}(\hat{\omega}) | N_s \leq A \text{ for all } s \text{ but one}\})$$

$$= 12^{-\hat{N}} \sum_{t=1}^{R} d_{k_t} \prod_{s \neq t} \sum_{N_s = 0}^{A} D_{N_s}^{(k_s)} 12^{-N_s}.$$

On the other hand, for fixed $1 \leq t, u \leq R, t \neq u$, we have

$$D_N^{-1} \sum_{\substack{N_1 + \cdots + N_R = N - \hat{N} \\ N_t \geq (N - \hat{N})/R, N_u \geq A}} \prod_{s=1}^{R} D_{N_s}^{(k_s)}$$

$$(25) \quad \leq \text{const} \cdot \sum_{\substack{N_s, s \neq t, u \\ N_u \geq A}} 12^{-\hat{N}} \left(\frac{NR}{N - \hat{N}}\right)^{5/2} N_u^{-5/2} \prod_{s \neq t, u} D_{N_s}^{(k_s)} 12^{-N_s}$$

$$\leq \text{const} \cdot A^{-3/2} \prod_{s \neq t, u} w_{k_s}$$

$$= \text{const} \cdot A^{-3/2},$$

where the constants depend on $\hat{\omega}$ only and we have used (7). This shows that

(26) $\quad \mu_N(\{\omega \in \mathcal{B}_{1/r}(\hat{\omega}) | \exists u \neq t \text{ s.t. } N_u \geq A, N_t \geq A\}) \leq \text{const} \cdot A^{-3/2}$,

where the constant depends on $\hat{\omega}$ only. Letting $A \to \infty$, we finally conclude from (23), (24) and (26) that

(27) $$\mu_N(\mathcal{B}_{1/r}(\hat{\omega})) \stackrel{N\to\infty}{\longrightarrow} 12^{-|\hat{\omega}|} \sum_{t=1}^{R} d_{k_t} \prod_{s \neq t} w_{k_s}.$$

This concludes the proof of Theorem 3.1, since it is clear from the definition that the (countable) set of finite well-labeled trees has vanishing $\mu$-measure. □

For later use we note that the proof extends immediately to the corresponding situation for $k$-labeled trees. Defining $\mu_N^{(k)}$ as the uniform probability measure concentrated on the set $\mathcal{W}'^{(k)}_N$ of $k$-labeled trees of size $N$, that is,

$$\mu_N^{(k)}(\omega) = (D_N^{(k)})^{-1} \qquad \text{for } \omega \in \mathcal{W}'^{(k)}_N, \qquad \mu_N^{(k)}(\overline{\mathcal{W}}^{(k)} \setminus \mathcal{W}'^{(k)}_N) = 0,$$

we thus have the corollary of the preceding proof.

COROLLARY 3.2. *The sequence $\mu_N^{(k)}$, $N \in \mathbb{N}$, converges to a Borel probability measure, $\mu^{(k)}$, concentrated on $\mathcal{W}_\infty^{(k)}$. We call $\mu^{(k)}$ the uniform probability measure on $\mathcal{W}_\infty^{(k)}$.*

Having proven the existence of the measure $\mu$, we easily get the following result on the measure $d\mu(\omega_1, \ldots, \omega_R | A(\hat{\omega}))$ that is obtained by conditioning $\mu$ on the event

$$A(\hat{\omega}) = \mathcal{B}_{1/r}(\hat{\omega}) = \{\omega \in \overline{\mathcal{W}} | B_r(\omega) = \hat{\omega}\},$$

where $\hat{\omega} \in \mathcal{W}_r$ is a finite tree of height $r$ and with $R$ vertices at maximal distance $r$ from the first root $i_0$, and we identify $A(\hat{\omega})$ (homeomorphically) with $\overline{\mathcal{W}}^{(k_1)} \times \cdots \times \overline{\mathcal{W}}^{(k_R)}$ as previously, where $k_1, \ldots, k_R$ are the labels of those $R$ vertices.

COROLLARY 3.3. *For $\hat{\omega} \in \mathcal{W}_r$, we have*

(28) $$\mu(A(\hat{\omega})) = 12^{-|\hat{\omega}|} \sum_{t=1}^{R} d_{k_t} \prod_{s \neq t} w_{k_s}$$

*and*

$$(29) \quad d\mu(\omega_1,\ldots,\omega_R|A(\hat{\omega})) = \mu(A(\hat{\omega}))^{-1} \sum_{t=1}^{R} d\mu^{(k_t)}(\omega_t) 12^{-|\hat{\omega}|} \prod_{s \neq t} d\rho^{(k_s)}(\omega_s),$$

*where the measure $\rho^{(k)}$ is supported on $\mathcal{W}^{(k)}$ and defined by*
$$\rho^{(k)}(\omega) = 12^{-|\omega|} \qquad \text{for } \omega \in \mathcal{W}^{(k)}.$$

PROOF. The identity (28) is an immediate consequence of (27).
To establish (29) it suffices, by Theorem 2.2 in [9], to prove

$$(30) \qquad \lim_{N \to \infty} \mu_N(\mathcal{B}_{1/r_1}(\hat{\omega}_1) \times \cdots \times \mathcal{B}_{1/r_R}(\hat{\omega}_R))$$

$$(31) \qquad = \sum_{t=1}^{R} \mu^{(k_t)}(\mathcal{B}_{1/r_t}(\hat{\omega}_t)) 12^{-|\hat{\omega}|} \prod_{s \neq t} \rho^{(k_s)}(\mathcal{B}_{1/r_s}(\hat{\omega}_s))$$

for any $r_1, \ldots, r_R \in \mathbb{N}$ and $\hat{\omega}_1, \ldots, \hat{\omega}_R \in \mathcal{W}$. Here $\mathcal{B}_{1/r_1}(\hat{\omega}_1) \times \cdots \times \mathcal{B}_{1/r_R}(\hat{\omega}_R)$ is considered as a subset of $A(\hat{\omega})$ by the previously mentioned homeomorphism. As above we may assume $\hat{\omega}_i \in \mathcal{W}_{r_i}, i = 1, \ldots, R$. Letting $\tilde{\omega}$ denote the tree obtained by grafting $\hat{\omega}_1, \ldots, \hat{\omega}_R$ at the leaves at maximal height in $\hat{\omega}$, then $\mathcal{B}_{1/r_1}(\hat{\omega}_1) \times \cdots \times \mathcal{B}_{1/r_R}(\hat{\omega}_R)$ consists of those trees in $\bar{\mathcal{W}}$ obtained by grafting arbitrary trees at the leaves of $\tilde{\omega}$, which correspond to leaves in some $\hat{\omega}_i$ at maximal height. Having made this observation, the proof of (31) is essentially identical to that of (27), taking into account the identity (28) and its analog for $\mu^{(k)}$; see (32) below. □

REMARK 3.4. For $\hat{\omega} \in \mathcal{W}_r^{(k)}$ and
$$A(\hat{\omega}) = \{\omega \in \overline{\mathcal{W}}^{(k)} | B_r(\omega) = \hat{\omega}\},$$
we have similarly with the same notation that

$$(32) \qquad \mu^{(k)}(A(\hat{\omega})) = 12^{-|\hat{\omega}|} \sum_{t=1}^{R} \frac{d_{k_t}}{d_k} \prod_{s \neq t} w_{k_s}$$

and

$$d\mu^{(k)}(\omega_1, \ldots, \omega_R | A(\hat{\omega})) = \mu^{(k)}(A(\hat{\omega}))^{-1} \sum_{t=1}^{R} d\mu^{(k_t)}(\omega_t) 12^{-|\hat{\omega}|} \prod_{s \neq t} d\rho^{(k_s)}(\omega_s).$$

REMARK 3.5. It is worth noting that the proof of the existence of the limit (27) and of the factorized form (29) of the conditional probability measure depends crucially on the fact that the exponent $-\frac{5}{2}$ of the asymptotic power dependence of $D_N^{(k)}$ on $N$ has value less than $-1$. For unlabeled planar trees this exponent is $-\frac{3}{2}$ such that similar, but simpler, arguments apply; see [15, 16].

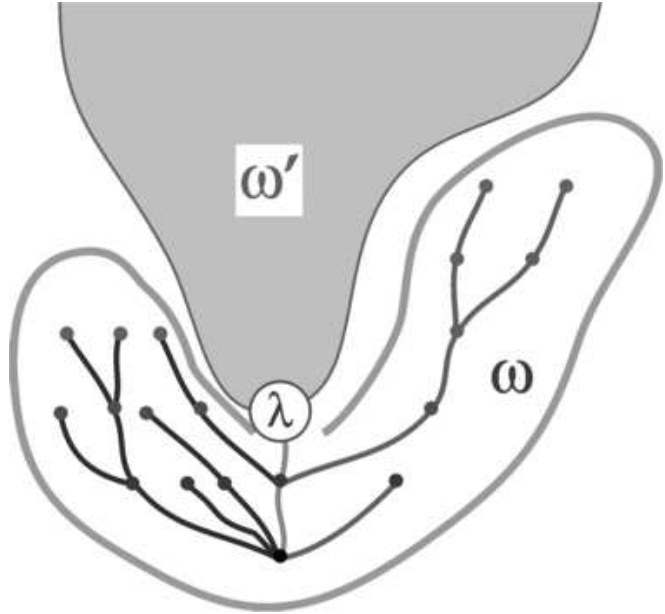

Fig. 2. *Grafting a tree:* $\mathrm{gr}(\omega, \lambda, \omega')$.

## 4. Description of the uniform probability measure on $\mathcal{W}_\infty$.

4.1. *More on the topology on labeled trees.* Consider a labeled tree $\omega \in \mathcal{W}^{(k)}$ with a marked leaf $\lambda$. Let $h(\lambda)$ denote the height of $\lambda$, that is, the graph distance in $\omega$ between $\lambda$ and the first root, and let $\ell(\lambda)$ be the label of $\lambda$. Let $\mathrm{gr}(\omega, \lambda, \omega')$ denote the labeled tree formed, starting from $\omega$, by grafting a labeled tree $\omega' \in \overline{\mathcal{W}}^{(\ell(\lambda))}$ at the leaf $\lambda$ of $\omega$ (see Figure 2). Set

$$\mathrm{Gr}(\omega, \lambda) = \{\mathrm{gr}(\omega, \lambda, \omega') | \omega' \in \overline{\mathcal{W}}^{(\ell(\lambda))}\}.$$

PROPOSITION 4.1. (i) *Each set* $\mathrm{Gr}(\omega, \lambda)$ *is both open and closed, and can be written as a union of sets in* $\mathcal{V}$.

(ii) *We have*

$$\mu(\mathrm{Gr}(\omega, \lambda)) = 12^{-|\omega|} d_{\ell(\lambda)}.$$

(iii) *If* $\hat{\omega}$ *is finite with radius* $r$, *then* $\mathcal{B}_{1/r}(\hat{\omega}) = A(\hat{\omega})$ *can be written, up to a set of $\mu$-measure* 0, *as a union of sets* $\mathrm{Gr}(\omega, \lambda)$, *all satisfying* $h(\lambda) = r$.

PROOF. (i) Setting

$$a = \rho(\omega) - h(\lambda) + 1 \quad \text{and} \quad k = \ell(\lambda),$$

we have
$$\mathrm{Gr}(\omega,\lambda) = \bigcup_{\substack{\omega' \in \mathcal{W}^{(k)} \\ \rho(\omega') \leq a}} \mathcal{B}_{1/(\rho(\omega)+1)}(\mathrm{gr}(\omega,\lambda,\omega')).$$

Since the family of all open balls with a given radius forms a partition of $\overline{\mathcal{W}}$, both $\mathrm{Gr}(\omega,\lambda)$ and its complement are open, being unions of open balls with radius $\frac{1}{\rho(\omega)+1}$.

(ii) By (i) we have
$$\lim_N \mu_N(\mathrm{Gr}(\omega,\lambda)) = \mu(\mathrm{Gr}(\omega,\lambda)).$$

Here
$$\mu_N(\mathrm{Gr}(\omega,\lambda)) = \frac{D^{(k)}_{N-|\omega|}}{D_N} = \frac{D^{(k)}_{N-|\omega|}}{D^{(k)}_N} \frac{D^{(k)}_N}{D_N}$$

for $N > |\omega|$ and (ii) follows from Corollary 2.6.

(iii) Given a tree $\hat{\omega}$, the grafting operation can be generalized in an obvious way to a finite sequence of leaves $\underline{\lambda} = (\lambda_1, \ldots, \lambda_R)$ in $\hat{\omega}$ and a finite sequence of trees $\underline{\omega} = (\omega_1, \ldots, \omega_R)$ with $\omega_i \in \overline{\mathcal{W}}^{(\ell(\lambda_i))}$. We denote the resulting tree by $\mathrm{gr}(\hat{\omega}, \underline{\lambda}, \underline{\omega})$.

Assume that the tree $\hat{\omega}$ has height $r$ and, as in Section 3, let $\Xi_r(\hat{\omega}) = ((j_1, k_1), \ldots, (j_R, k_R))$. Furthermore, let $\Lambda_t$ be the sequence obtained by erasing $j_t$ from $(j_1, \ldots, j_R)$. Setting
$$\overline{\Omega}_t = \prod_{s \neq t} \overline{\mathcal{W}}^{(k_s)}, \qquad \Omega_t = \prod_{s \neq t} \mathcal{W}^{(k_s)},$$

we then have
$$A(\hat{\omega}) = \bigcup_{t=1}^R \bigcup_{\underline{\omega} \in \overline{\Omega}_t} \mathrm{Gr}(\mathrm{gr}(\hat{\omega}, \Lambda_t, \underline{\omega}), j_t)$$
$$\supset \bigcup_{t=1}^R \bigcup_{\underline{\omega} \in \Omega_t} \mathrm{Gr}(\mathrm{gr}(\hat{\omega}, \Lambda_t, \underline{\omega}), j_t) \doteq A'(\hat{\omega}).$$

Thus, the set $A'(\hat{\omega})$ consists of the trees in $A(\hat{\omega})$ such that only one among the $R$ grafted trees is infinite, the others being finite. From (ii) and (28) we conclude that
$$\mu(A'(\hat{\omega})) = \sum_{t=1}^R \sum_{(\omega_1, \ldots, \omega_{R-1}) \in \Omega_t} d_{k_t} 12^{-|\hat{\omega}|-|\omega_1|-\cdots-|\omega_{R-1}|}$$
$$= \mu(A(\hat{\omega})). \qquad \square$$

REMARK 4.2. Note that statement (ii) in Proposition 4.1 has a straightforward generalization to a family of leaves, possibly at different heights.

4.2. *The spine.* We define a *spine* of a labeled tree $\omega$ to be any infinite sequence of labeled vertices $(\lambda_r)$ that start at the first root vertex of $\omega$ and such that $\lambda_r$ is a son of $\lambda_{r-1}$. Equivalently, a spine of $\omega$ is an infinite labeled linear subtree with the same first root vertex as $\omega$. Similarly, for $r \in \mathbb{N}$, an $r$-spine of $\omega$ is a labeled linear subtree of height $r$ with same first root vertex as $\omega$.

The following result will be important for the subsequent developments.

THEOREM 4.3. *With $\mu$-probability* 1, *a tree contains exactly one spine.*

In other words, the probability measure $\mu$ is supported on the set $\mathcal{S}$ of trees with exactly one spine. Such trees are obtained by grafting, at each vertex of its spine, a pair of *finite* well-labeled trees, one on the right and one on the left. This is not unexpected, since the limit law in the simpler case of unlabeled trees has a similar description [1, 15, 16, 19, 24].

Note that by convention we draw the root at the bottom of the tree so the right of the spine means the right, when looking at the spine *from the root*.

PROOF OF THEOREM 4.3. Set

$$\mathcal{S}_r = \bigcup_{\substack{\omega \in \mathcal{W} \\ h(\lambda)=r}} \mathrm{Gr}(\omega,\lambda) \quad \text{and} \quad \mathcal{S}' = \bigcap_{r \geq 1} \mathcal{S}_r.$$

Due to Proposition 4.1(iii), $\mu(\mathcal{S}_r) = 1$ and thus $\mu(\mathcal{S}') = 1$.

Assume now that $\eta$ belongs to $\mathcal{S}'$. Then $\eta$ has infinite size, because it contains a vertex at each height $r$. Being an element in $\mathcal{S}_r \cap \mathcal{W}_\infty$, $\eta$ has a unique decomposition of the form $\mathrm{gr}(\omega_r, \lambda_r, \omega'_r)$, where $\omega_r \in \mathcal{W}$ and $h(\lambda_r) = r$. Necessarily, in such a decomposition, $\omega'_r$ belongs to $\mathcal{W}_\infty^{(\ell(\lambda_r))}$. Hence, for each $r$, there is a unique pair $(\omega_r, \lambda_r)$ such that $\eta \in \mathrm{Gr}(\omega_r, \lambda_r)$. It follows that all spines in $\eta$, if any, coincide up to height $r$ and that $\lambda_{r+1}$ is necessarily the son of $\lambda_r$. This shows that $(\lambda_r)_{r \geq 0}$ is the unique spine of $\eta$, so that $\mathcal{S}' \subset \mathcal{S}$. □

For a random $\mu$-distributed element $\omega$ of $\mathcal{S}(\subset \mathcal{W}_\infty)$, we introduce the following notation. By $e_n$ we denote the vertex at height $n$ on its spine, and the label of $e_n$ is denoted by $X_n(\omega)$. Furthermore, we let $L_n(\omega)$ [resp. $R_n(\omega)$] be the *finite* subtree of $\omega$ attached to $e_n$ on the left (resp. on the right) of its spine.

By $\hat{\rho}^{(k)}$ we denote the measure obtained by normalizing $\rho^{(k)}$, that is,

$$\hat{\rho}^{(k)}(\omega) \doteq \frac{12^{-|\omega|}}{w_k} \quad \text{for } \omega \in \mathcal{W}^{(k)}.$$

Also, for any sequence of triples of nonnegative numbers $\{(p_k, r_k, q_k) : k \geq 1\}$ such that $p_k + r_k + q_k = 1$ and $q_1 = 0$, any Markov chain $\{Z_k : k \geq 0\}$ on the positive integers, with transition probabilities defined by

$$\mathbb{P}(Z_{n+1} = k+1 | Z_n = k) = p_k,$$
$$\mathbb{P}(Z_{n+1} = k-1 | Z_n = k) = q_k,$$
$$\mathbb{P}(Z_{n+1} = k | Z_n = k) = r_k,$$

will be called a *birth and death process* with parameters $\{(p_k, r_k, q_k) : k \geq 1\}$. (Usually such a birth and death process is defined on the nonnegative integers [20] and one assumes that $q_0 = 0$.)

THEOREM 4.4. *The measure $\mu$ has the following probabilistic description.*

(i) *The variable $X = (X_n)_{n \geq 0}$ is a birth and death process with parameters $\{(p_k, r_k, q_k) : k \geq 1\}$ defined by*

$$(33) \qquad q_k \doteq \frac{(w_k)^2}{12 d_k} d_{k-1}, \qquad r_k \doteq \frac{(w_k)^2}{12}, \qquad p_k \doteq \frac{(w_k)^2}{12 d_k} d_{k+1},$$

*where by convention $d_0 = 0$.*

(ii) *Conditionally, given that $X = (s_n)_{n \geq 0}$, the $L_n$'s and $R_n$'s form two independent sequences of independent random labeled trees, distributed according to the measures $\hat{\rho}^{(s_n)}$.*

PROOF. Given $r \in \mathbb{N}_0$ and $(s_0, s_1, \ldots, s_r) \in \mathbb{N}^r$ fulfilling

$$s_0 = 1 \quad \text{and} \quad |s_k - s_{k-1}| \leq 1 \qquad \text{for } k = 1, \ldots, r,$$

we observe that the set

$$\{\omega \in \mathcal{S} | (L_k, R_k) = (\omega'_k, \omega''_k), 0 \leq k \leq r-1, \text{ and } X_k = s_k, 0 \leq k \leq r\}$$

equals $\mathrm{Gr}(\omega, \lambda) \cap \mathcal{S}$, where $\omega$ is formed by grafting the pairs $(\omega'_k, \omega''_k)$ of *finite* well-labeled trees in $\mathcal{W}^{(s_k)}$ on both sides of the $r$-spine $((e_0, s_0), (e_1, s_1), \ldots, (e_r, s_r))$ and where $\lambda = (e_r, s_r)$.

Thus, by Proposition 4.1(ii), the "finite-dimensional distributions" are given by

$$\mu(\{(L_k, R_k, X_k) = (\omega'_k, \omega''_k, s_k), 0 \leq k \leq r-1, \text{ and } X_r = s_r\})$$
$$= d_{s_r} 12^{-|\omega|}$$
$$= \prod_{k=0}^{r-1} \frac{1}{12} \frac{d_{s_{k+1}}}{d_{s_k}} (w_{s_k})^2 \prod_{k=0}^{r-1} \frac{12^{-|\omega'_k|}}{w_{s_k}} \frac{12^{-|\omega''_k|}}{w_{s_k}}.$$

Upon realizing that relationships (15)–(17) imply $p_k + q_k + r_k = 1$, the statements in the theorem can be read off from this formula. $\square$

REMARK 4.5. By a slight extension of the argument, it follows that, upon conditioning on a fixed $r$-spine, all $2r + 1$ branches, including the infinite one attached to the end of the $r$-spine, are independently distributed, and the latter is distributed according to the measure $\mu_{s_r}$, where $s_r$ is the end label (cf. also Corollary 3.3).

4.3. *The branches.* The following theorem describes in more detail the probabilistic structure of finite subtrees grafted on the left and on the right of the spine: the labels of the nodes can be seen as types of *multitype* Galton–Watson processes. In this setting we have the following theorem:

THEOREM 4.6. *Conditionally, given that $X_n = k$, $R_n$ and $L_n$ are independent multitype Galton–Watson trees, in which the ancestor has type $k$ and an (type $\ell$) individual can only have progeny of type $\ell + \varepsilon$, $\varepsilon \in \{0, \pm 1\}$. In such multitype Galton–Watson trees, the progeny of an (type $\ell$) individual is determined by a sequence of independent trials with four possible outcomes, $\ell + 1$, $\ell - 1$, $\ell$ and $e$ (for "extinction"), with respective probabilities, $w_{\ell+1}/12$, $w_{\ell-1}/12$, $w_\ell/12$ and $\frac{1}{w_\ell}$ (that add up to 1), that are sequence-stopped just before the first occurrence of $e$: an (type $\ell$) individual has as many children of type $\ell + 1$, $\ell - 1$ or $\ell$ as there are occurrences of $\ell + 1$, $\ell - 1$ and $\ell$ in the sequence before the first occurrence of $e$.*

PROOF. Note that, owing to (18),

$$(34) \qquad (w_{\ell+1} + w_{\ell-1} + w_\ell)\frac{1}{12} + \frac{1}{w_\ell} = 1.$$

With respect to $\hat{\rho}^{(k)}$, the probability that the ancestor has $m$ sons with respective labels $(k_1, \ldots, k_m)$ and with associated subtrees $(\omega_1, \ldots, \omega_m)$ is

$$\frac{12^{-m-\sum |\omega_i|}}{w_k},$$

provided that $|k - k_i| \leq 1$ and $\omega_i \in \mathcal{W}_{k_i}$. This probability can be written

$$\frac{1}{w_k} \prod_{i=1}^m \frac{w_{k_i}}{12} \prod_{i=1}^m \hat{\rho}^{(k_i)}(\omega_i).$$

The theorem follows by induction. □

REMARK 4.7. Thus, for the special case of multitype Galton–Watson (GW) trees we consider here, the picture is quite like the picture given in [1, 15, 16, 19, 22, 24] for monotype Galton–Watson trees. However, there are some differences: the spine is not pasted uniformly on the available leaves, but rather with a bias introduced by the different types. Also, the finite branching trees are multitype GW trees, critical in the sense that the average progeny of an (type $\ell$) individual,

$$w_\ell - 1 = 1 - \frac{4}{(\ell+1)(\ell+2)},$$

has supremum in $\ell$ equal to 1, but we can see from the expansion of $W^{(\ell)}(z)$ at $1/12$ that, given that the label of its root is $\ell$, the conditional expected size of such a "critical" GW tree is finite [and equal to $(3\ell^2 + 9\ell - 2)/10$], while the (finite) size of a standard critical GW tree usually has an infinite expectation. As an additional feature, the succession of types on the spine is a birth and death process. To be complete, in [22], a probability measure on infinite GW trees with a finite number of types, quite similar to $\mu$, is used to give an elegant proof of the Kesten–Stigum theorem.

**5. Label occurrences in a uniform labeled tree.** As will be seen in Section 6 the volume of the ball with radius $k$ in a random uniform quadrangulation with $N$ faces has the same distribution as the number of nodes with label *smaller than* $k+1$ in a well-labeled tree with $N$ edges. In this section we study the number of nodes with label *exactly* $k$ in the uniform infinite well-labeled tree, that is, the number $N_k$ of occurrences of label $k$ in a well-labeled tree with respect to the measure $\mu$. In particular, we determine the asymptotic behavior of the average value $\mathbb{E}_\mu[N_k]$. For this purpose we need to investigate, in Sections 5.1 and 5.2, two types of random walks associated to $\mu$.

5.1. *The random walk along the spine.* We wish to determine the asymptotic behavior of the number $S_k$ of occurrences of label $k$ along the spine of the uniform infinite well-labeled tree. This behavior eventually depends on the asymptotic behavior of $q_k, r_k$ and $p_k$ as $k \to \infty$. We shall prove the next proposition:

PROPOSITION 5.1. *We have*
$$\lim_{k\to\infty} \frac{\mathbb{E}_1[S_k]}{k} = \frac{3}{7}.$$

In probabilistic terms, $S_k$ can be seen as the sojourn time of the process $X$ at level $k$, that is,

$$S_k(\omega) = \sharp\{n \in \mathbb{N} | X_n(\omega) = k\}.$$

Proposition 5.1 will be obtained as a consequence of a more general property of sojourn times $S_k$ of general birth and death processes, stated in the next lemma.

LEMMA 5.2. *Assume:*

(i) *All $q_k$ and $p_k$ are positive with the exception of $q_1 = 0$.*
(ii) *There exist constants $a > 1$ and $c$ such that, for all $k$,*

$$k^2 \left| \frac{q_k}{p_k} - 1 + \frac{a}{k} \right| \leq c.$$

(iii) *The sequence $r_k$, $k \in \mathbb{N}$, has a limit $b < 1$.*

*Then*

$$\lim_k \frac{\mathbb{E}_i[S_k]}{k} = \frac{2}{(1-b)(a-1)}.$$

REMARK 5.3. In (ii), the assumption $a > 1$ ensures that $X$ is transient and that a.s. $\lim X_n = +\infty$.

REMARK 5.4. As is usual for Markov chains, the index $i$ in $\mathbb{E}_i[S_k]$ stresses the assumption $\mathbb{P}(X_0 = i) = 1$. For instance, in the case of the random uniform infinite labeled tree, $\mathbb{E}_i[S_k] = \langle S_k \rangle_{\mu^{(i)}}$. Due to transience, $\mathbb{E}_i[S_k]$ does not depend on $i$ if $i \leq k$.

PROOF OF LEMMA 5.2. Set

$$m_{i,k} = \frac{q_i q_{i+1} \cdots q_{i+k-1}}{p_i p_{i+1} \cdots p_{i+k-1}},$$

which is positive for $i \geq 2$ by assumption. By definition, set $m_{i,0} = 1$ and let $T_i$ denote the first hitting time of level $i$. As is well known ([20], Chapter 3), $X$ is transient if and only if $m_{1,k}$ is the general term of a converging series. Also

$$\mathbb{P}_{i+k}(T_i = +\infty) = \frac{\sum_{j=0}^{k-1} m_{i+1,j}}{\sum_{j \geq 0} m_{i+1,j}},$$

with the special case

$$\mathbb{P}_{i+1}(T_i = +\infty) = \frac{1}{\sum_{j \geq 0} m_{i+1,j}} \doteq \rho_i.$$

Let $D_k$ denote the number of downcrossings $k+1 \downarrow k$ and let $Y_\ell$ be the sojourn time of $X$ at level $k$ after the $\ell$th and before the $(\ell+1)$st (if it occurs) downcrossing. Also, let $Y_0$ be the sojourn time of $X$ at level $k$ before the first (if any) downcrossing. Then, for $i \leq k$ and $n \geq 0$,

$$\mathbb{P}_i(D_k \geq n) = (1 - \rho_k)^n.$$

Furthermore, conditionally, given that $D_k = n$, $(Y_0, Y_1, \ldots, Y_n)$ is a sequence of i.i.d. random variables that satisfy, for $m \geq 1$,

$$\mathbb{P}_i(Y_j \geq m) = (1 - p_k)^{m-1}.$$

As a consequence, by Wald's identity,

$$\mathbb{E}_i[S_k] = \mathbb{E}_i[1 + D_k]\mathbb{E}_i[Y_0] = \frac{1}{\rho_k p_k}. \tag{35}$$

To determine the asymptotics of $\rho_k$, we use that, due to assumption (ii), the expression

$$k^2\left(\frac{q_k(k+1)^a}{p_k k^a} - 1\right)$$

is bounded as a function of $k$. Thus

$$\limsup_k \sup_{\ell \geq 0}\left|m_{k,\ell}\left(\frac{k+\ell}{k}\right)^a - 1\right| = 0$$

and, as a consequence,

$$\lim_k \frac{1}{k}\sum_{\ell \geq 0} m_{k,\ell} = \lim_k \frac{1}{k}\sum_{\ell \geq 0} \frac{1}{(1+\ell/k)^a} = \int_1^{+\infty} \frac{dx}{x^a} = \frac{1}{a-1}.$$

Combining this with the definition of $\rho_k$ and

$$\lim_k p_k = \frac{1-b}{2},$$

the claimed limit follows from (35). □

In the next lemma we collect the large-$k$ behaviors of $q_k, r_k$ and $p_k$ that follow immediately from (10), (14) and (33).

LEMMA 5.5. *We have*

$$q_k = \frac{1}{3} - \frac{4}{3k} + \mathcal{O}(k^{-2}),$$

$$r_k = \frac{1}{3}\left(1 - \frac{4}{k^2}\right) + \mathcal{O}(k^{-3}),$$

$$p_k = \frac{1}{3} + \frac{4}{3k} + \mathcal{O}(k^{-2}).$$

Proposition 5.1 now follows from Lemma 5.2 and Lemma 5.5, with $a = 8$ and $b = 1/3$.

REMARK 5.6. By discretization of the $d$-dimensional Bessel process, one obtains a birth and death process that satisfies the assumptions of Lemma 5.2 for $(a,b) = (d-1, 0)$. So the birth and death process we meet here is, in a sense, close to a nine-dimensional Bessel process. Similar birth and death processes also appear, in connection with random nonlabeled trees, in the study of the three-dimensional Bessel process and lead to an elegant proof of the decomposition theorem of Williams [23].

Proposition 5.1 ensures that $S_k$ is finite almost surely; in other words, that each label $k$ occurs almost surely only finitely many times on the spine. It will turn out to be essential for the interpretation of well-labeled trees as quadrangulated surfaces in the next section that the same result holds for the total number $N_k$ of occurrences of a label $k$ on the whole tree, that is, that $N_k$ is finite almost surely for each $k \in \mathbb{N}$. In Section 5.3, we shall prove that $\mathbb{E}_\mu[N_k]$ is finite and we shall also determine its asymptotic behavior as $k \to \infty$. For this purpose, a study of the average value of $N_k$ in a generic branch attached to the spine is first needed.

5.2. *A random walk associated to branches.* Consider a branch $\omega_n$, left or right, grafted at the $n$th site $e_n$ of the spine of a well-labeled tree $\omega \in \mathcal{S}$. According to Theorem 4.4, conditionally, given that $k$ is the label of $e_n$, $\omega_n$ is distributed according to $\hat{\rho}^{(k)}$. Letting $N_j(\omega_n)$ denote the number of occurrences of the label $j$ in $\omega_n$, we denote by $G(k,j)$ the (normalized) average value of $N_j$, that is,

$$G(k,j) = \mathbb{E}_{\hat{\rho}^{(k)}}[N_j]$$
$$= \frac{1}{w_k} \sum_{\omega \in \mathcal{W}^{(k)}} N_j(\omega) 12^{-|\omega|}.$$

Consider the set $\mathcal{W}_\star^{(k)}$ of *marked* finite trees in $\mathcal{W}^{(k)}$, that is, couples $(\omega, e)$ with $\omega \in \mathcal{W}^{(k)}$ and $e$ a marked vertex of $\omega$, endowed with the measure $\frac{12^{-|\omega|}}{w_k}$ for each element $(\omega, e)$ of $\mathcal{W}_\star^{(k)}$. Then $G(k,j)$ can be seen as the measure of the set of marked trees whose marked vertex has label $j$.

Given an element $(\omega, e)$ of $\mathcal{W}_\star^{(k)}$, it has a distinguished (finite) spine

$$((f_0, \theta_0), (f_1, \theta_1), \ldots, (f_L, \theta_L)),$$

namely the path connecting the root $i_0(\omega) = f_0$ to the marked vertex $e = f_L$. We let

$$(\omega'_1, \ldots, \omega'_{L-1}) \quad \text{and} \quad (\omega''_1, \ldots, \omega''_{L-1}),$$

respectively, with $\omega'_t$ and $\omega''_t$ in $\mathcal{W}^{(\theta_t)}$, denote the sequences of subtrees (branches) in $\omega$ grafted on this spine on the left and right, and denote by $\omega_L \in \mathcal{W}^{(j)}$ the subtree attached to the marked vertex. Since this correspondence between marked trees on one side and spines together with branches on the other is bijective, we obtain, by summing over branches first, the following representation:

LEMMA 5.7. *We have*

(36) $$G(k,j) = \sum_{\theta \,:\, k \to j} 12^{-|\theta|} \prod_{t=0}^{|\theta|-1} w_{\theta_t} w_{\theta_{t+1}}$$

*in which*

$$\theta = (\theta_0, \theta_1, \ldots, \theta_L) \quad \text{and} \quad |\theta| = L$$

*and the sum is over spine-label sequences $\theta$ with initial label $k$ and final label $j$.*

Evidently, $\theta$ can also be viewed as a walk in $\mathbb{N}$ from $k$ to $j$ whose steps $\theta_{t+1} - \theta_t$ belong to $\{0, \pm 1\}$ and have "probability"

$$\frac{w_{\theta_t} w_{\theta_{t+1}}}{12}.$$

As an alternative proof of (36), by Theorem 4.6, the average number of children of an individual of type $k$ is $w_k - 1$, so that for $\varepsilon \in \{0, \pm 1\}$, due to Wald's identity, the average number $p_{k,k+\varepsilon}$ of children of type $k + \varepsilon$ of an individual of type $k$ is given by

$$p_{k,k+\varepsilon} = (w_k - 1) \frac{w_{k+\varepsilon}}{w_{k-1} + w_k + w_{k+1}}$$

$$= \frac{w_k w_{k+\varepsilon}}{12},$$

the second equality owing to (34). Let $\Sigma_\theta(\omega_n)$ denote the number of spines of $\omega_n$ whose spine-label sequence is $\theta$. Let us consider $\theta$ as a word in

$$\mathbb{N}^\star = \{\varnothing\} \cup \left(\bigcup_{k \geq 1} \mathbb{N}^k\right).$$

Then, for any $a$ in $\mathbb{N}$, conditionally given $\Sigma_\theta = k$, $\Sigma_{\theta a}$ is distributed as

$$\sum_{i=1}^k Y_i,$$

in which $Y_i$ stands for the number of children of type $a$ of the last vertex of the $i$th spine with spine-label sequence $\theta$. Thus the $Y_i$'s are i.i.d. with expectation $p_{\theta_L, a}$ and, according to Wald's identity,

$$\mathbb{E}[\Sigma_{\theta a}] = \mathbb{E}[\Sigma_\theta]\mathbb{E}[Y_1] = \mathbb{E}[\Sigma_\theta] p_{\theta_L, a}.$$

By induction, the average number of paths $\theta$ in a tree distributed according to $\hat{\rho}^{(k)}$ is given by

$$\prod_{t=0}^{L-1} p_{\theta_t, \theta_{t+1}} = 12^{-L} \prod_{t=0}^{L-1} w_{\theta_t} w_{\theta_{t+1}}.$$

Summing over all spines with initial label $k$ and final label $j$, we obtain again the formula for $G(k, j)$.

THEOREM 5.8. *The function G has the following properties:*

(i) $G(k,j) = G(j,k)$ *for* $k,j \in \mathbb{N}$,
(ii) $\phi'_+(k)\phi'_-(j) \leq G(k,j) \leq \phi_+(k)\phi_-(j)$ *for* $j > k \geq 1$,
(iii) $G(k,k) = \Theta(k)$,

*where*

(37) $$\phi'_+(k) = \Theta(k^4), \qquad \phi_+(k) = \Theta(k^4),$$

(38) $$\phi'_-(k) = \Theta(k^{-3}), \qquad \phi_-(k) = \Theta(k^{-3}).$$

PROOF. The symmetry property (i) follows from (36) (and from $p_{i,j} = p_{j,i}$). From (36) we also deduce the difference equation

$$G(k,j) = \delta_{kj} + \tfrac{1}{12}w_j \sum_{\varepsilon=0,\pm 1} w_{j+\varepsilon} G(k, j+\varepsilon),$$

where we use the convention $w_0 = 0$. Setting

(39) $$H(k,j) = \tfrac{1}{12} w_k G(k,j) w_j,$$

this equation can be rewritten as

$$\Delta_j H(k,j) = 3(4(w_j)^{-2} - 1)H(k,j) - \delta_{kj},$$

where $\Delta_j$ is the discrete Laplace operator with respect to $j$, whose action on functions $\psi : \mathbb{N} \to \mathbb{R}$ is given by

$$\Delta_j \psi(j) = \psi(j+1) + \psi(j-1) - 2\psi(j),$$

in which, by convention, $\psi(0) = 0$.

We now proceed to establish the upper bound in (ii). Since $w_k \to 2$ as $k \to \infty$, we may replace $G$ by $H$. For any positive sequence $u_k$, $k \in \mathbb{N}$, let us define $H^u(k,j)$ analogously to $H$ by replacing $w$ by $u$ in formulas (36) and (39), so that $H^w = H$. Using

$$3(4w_k^{-2} - 1) \geq \frac{12}{k(k+3)}$$

and defining $u_k > 0$ by

$$3(4u_k^{-2} - 1) = \frac{12}{k(k+3)},$$

we have

(40) $$w_k \leq u_k < 2.$$

It follows that

$$H(k,j) \leq H^u(k,j)$$

and $H^u$ fulfills

(41) $$\Delta_j H^u(k,j) = \frac{12}{j(j+3)} H^u(k,j) - \delta_{kj}.$$

Two linearly independent solutions to

$$\Delta_j \psi(j) = \frac{12}{j(j+3)} \psi(j), \qquad j \geq 2,$$

are

$$\psi_+(j) = (j+3)(j+2)(j+1)j$$

and

$$\psi_-(j) = (j+3)(j+2)(j+1)j \sum_{k=j}^{\infty} \frac{1}{k(k+1)^2(k+2)^2(k+3)^2(k+4)},$$

leading to

(42) $$\psi_+(j) = \Theta(j^4) \quad \text{and} \quad \psi_-(j) = \Theta(j^{-3}).$$

We conclude that

$$H^u(k,j) = c_k \psi_-(j) + c_k^+ \psi_+(j) \qquad \text{for } j > k \geq 1,$$

where $c_k$ and $c_k^+$ are constants that depend only on $k$. Here $c_k^+ = 0$, because, according to (40), $H^u(k,j)$ is bounded above by the function $H^2(k,j)$, defined analogously to $H^u(k,j)$ by replacing $u_k$ by 2. This case corresponds to the simple random walk where $q_k = r_k = p_k = \frac{1}{3}$ for $k \geq 2$ with reflecting boundary condition at $k = 1$, and the function $H^2$ fulfills

$$\Delta_j H^2(k,j) = 0,$$

that is, it is a linear function of $j$. However, since $\psi_+(j) = \Theta(j^4)$, it follows that $c_k^+ = 0$ and so

$$H^u(k,j) = c_k \psi_-(j) \qquad \text{for } j > k \geq 1.$$

Using the symmetry of $H^u(k,j)$, we get that $c_k$ is a linear combination of $\psi_+(k)$ and $\psi_-(k)$ for $k \geq 2$, say

$$c_k = c \cdot \psi_-(k) + d \cdot \psi_+(k).$$

Since $H^u > 0$, we obviously have $d \geq 0$. We claim that $d > 0$. Otherwise, we would have

$$H^u(k,j) = c \cdot \psi_-(k) \psi_-(j) \qquad \text{for } k,j \geq 2, k \neq j$$

and, in particular, $H^u(k, k \pm 1) \to 0$ as $k \to \infty$. Using this in (41) for $k = j$, that is,

$$\left(2 + \frac{12}{k(k+3)}\right) H^u(k,k) = 1 + H^u(k,k-1) + H^u(k,k+1),$$

we conclude that $H^u(k,k) \to \frac{1}{2}$ as $k \to \infty$. However, this contradicts the inequality

$$H^u(k, k+1) \geq (u_{k+1})^2 H^u(k,k)/12,$$

which is an easy consequence of the definition of $H^u$. Hence we must have $d > 0$ and (42) implies

(43) $$H^u(k,k) = \Theta(k).$$

Defining $\phi_+(k) = c_k$ and $\phi_-(k) = \psi_-(k)$ for $k \in \mathbb{N}$, we have established the upper bound in (ii).

Similarly, using

$$3(4(w_j)^{-2} - 1) \leq \frac{12}{(k-1)(k+2)}, \qquad k \geq 2,$$

and defining $v_k \leq w_k$ by

$$v_k = \begin{cases} w_1, & \text{for } k = 1, \\ u_{k-1}, & \text{for } k \geq 2. \end{cases}$$

we obtain that $H^v$ is a lower bound for $H$ with the same asymptotic behavior as $H^u$. This finishes the proof of (ii). Finally, (iii) follows from (43) and the corresponding relationship for $H^v$. □

5.3. *Label occurrences.* We are now ready to compute the asymptotic behavior of the average number $\mathbb{E}_\mu[N_j]$ of occurrences of the label $j$ in the full uniformly distributed well-labeled tree.

THEOREM 5.9. *We have*

$$\mathbb{E}_\mu[N_j] = \Theta(j^3).$$

PROOF. For a random uniform well-labeled infinite tree $\omega \in \mathcal{S}$, we have

$$N_j(\omega) = \sum_{n=0}^{\infty}(N_j(R_n(\omega)) + N_j(L_n(\omega))) - S_j(\omega).$$

By Corollary 3.3 or Theorem 4.4(ii) we have

$$\sum_{n=0}^{\infty} \mathbb{E}_\mu[N_j(R_n) + N_j(L_n)] = 2\mathbb{E}_1\left[\sum_{n=0}^{\infty} G(X_n, j)\right]$$

$$= 2\mathbb{E}_1\left[\sum_{k=1}^{\infty} S_k G(k, j)\right]$$

so that

$$\mathbb{E}_\mu[N_j] = 2 \sum_{k=1}^{\infty} \mathbb{E}_1[S_k] G(k,j) - \mathbb{E}_1[S_j]$$

$$\leq 2 \sum_{k=1}^{j-1} \mathbb{E}_1[S_k] \phi_+(k) \phi_-(j)$$

$$+ 2 \sum_{k=j+1}^{\infty} \mathbb{E}_1[S_k] \phi_-(k) \phi_+(j)$$

$$+ \mathbb{E}_1[S_j](2G(j,j) - 1).$$

By combining with the asymptotic behaviors of $\phi_\pm$ and $\mathbb{E}_1[S_k]$ from (37), (38) and Proposition 5.1, one finds that the first two terms in the last expression equal $\Theta(j^3)$, while the last term equals $\Theta(j^2)$. Similarly, one obtains a lower bound of the same type, thus proving the theorem. □

COROLLARY 5.10. *For each $j \in \mathbb{N}$, the number of vertices with label $j$ is $\mu$-almost surely finite.*

**6. The uniform infinite random quadrangulation.** In this section we show how to draw an infinite planar map, starting from an infinite labeled tree in the set

$$\mathcal{C} = \{\omega \in \mathcal{S} | \forall j \geq 1, N_j(\omega) < +\infty\}$$

which, as we saw in the previous section, has $\mu$-measure 1. The result will be an infinite rooted quadrangulation of a domain in the plane, a notion defined more precisely below. Our construction follows the same steps as the finite analog given in [12], Section 3.4. In particular, the vertices of the tree can be identified with the vertices of the corresponding quadrangulation, with the exception of a distinguished root vertex in the latter. We stress, however, that we do not establish a bijection between infinite well-labeled trees and infinite quadrangulations in general. Actually, we think that the construction of a finite well-labeled tree from a finite quadrangulation extends to the infinite case, but we do not know if this last map is onto. What is important for our purposes, is that the constructed mapping $\mathcal{Q}$ allows us to transport the measure $\mu$ on $\mathcal{C}$ to the image set $\mathcal{Q}(\mathcal{C})$ of quadrangulations and, furthermore, that $\mathcal{Q}$ possesses the property that the label of a vertex of a tree $\omega \in \mathcal{C}$ equals the distance to the root of the corresponding vertex in the quadrangulation $\mathcal{Q}(\omega)$ (see Property 6.3 below).

6.1. *Regular infinite planar maps.* By an *infinite planar map* we mean an embedding $\mathcal{M} = E(G)$ of an infinite graph $G$, which we assume is connected and all of whose vertices are of finite degree, into the 2-sphere $S^2$, such that the edges are represented by smooth curve segments that do not intersect each other except at common vertices. To be a useful concept, some regularity properties of the embedding are necessary, in addition. For instance, it is evidently possible to embed the infinite linear tree into $S^2$, thought of as the plane $\mathbb{R}^2$ with the point $\infty$ added, such that it is mapped onto a circle, which we shall not accept as a valid embedding. To avoid this we make the following assumption (see also [8]):

ASSUMPTION $(\alpha)$. If $p_i$, $i \in \mathbb{N}$, is a sequence of points in a planar map $\mathcal{M}$, considered as the union of its edges in $S^2$, such that, for $i \neq j$, $p_i$ and $p_j$ are contained in different edges, then the sequence has no condensation point in $\mathcal{M}$.

Consider now a closed continuous curve $C$ in an infinite planar map $\mathcal{M}$ composed of a sequence of edges. By Assumption $(\alpha)$, the set of different edges must be finite and so the complement of $C$ in $\mathbb{R}^2$ decomposes into a finite number of connected components. If one of the components contains only a finite number of vertices of $\mathcal{M}$, then the part of $\mathcal{M}$ inside or on the boundary of this component is a finite planar map. The faces of this finite planar map inside the connected component are then also called faces of $\mathcal{M}$. Thus the faces of $\mathcal{M}$ are those obtained in this way for some closed curve $C$. In particular, each face is bounded by a polygonal loop composed of a *finite* number of edges. This leads us to make another assumption:

ASSUMPTION $(\beta)$. Either an edge of a planar map $\mathcal{M}$ is shared by exactly two faces or it occurs twice in the boundary of one face of $\mathcal{M}$.

For instance, such a regularity assumption does not hold for infinite trees, which have only one "face" with infinite degree. Let $\mathcal{D}(\mathcal{M})$ denote the union of the closed faces of $\mathcal{M}$: due to Assumption $(\beta)$, $\mathcal{D}(\mathcal{M})$ is an open connected subset of $S^2$. We identify embeddings $\mathcal{M}$ and $\mathcal{M}'$ that are related by an orientation-preserving homeomorphism between $\mathcal{D}(\mathcal{M})$ and $\mathcal{D}(\mathcal{M}')$.

The planar map is called a *quadrangulation* if all of its faces are quadrangles, that is, bounded by polygons with four edges. By a *rooted planar map* we mean a (finite or infinite) planar map with a distinguished oriented edge $(i_0, i_1)$, called the *root* of the planar map. As in the case of trees, we call $i_0$ the first root vertex.

6.2. *The mapping $\mathcal{Q}$.* To define $\mathcal{Q}$, we need regularity and uniqueness conditions on the embedding of labeled trees in $S^2$: for instance, even with the restriction of Assumption $(\alpha)$, the linear tree can be embedded in many nonhomeomorphic ways into $S^2$. However, requiring the embedding to be such that sequences of the type in Assumption $(\alpha)$ have exactly one and the same condensation point in $S^2$, then the embedding is unique up to homeomorphisms of $S^2$. More generally, this also holds for arbitrary infinite trees:

The combinatorial definition of a tree given at the beginning of Section 2 determines a unique embedding of the corresponding graph into $S^2$, up to homeomorphisms, such that all sequences as in Assumption $(\alpha)$ have exactly one and the same condensation point.

Indeed, such an embedding has already been indicated in Section 2, where vertices at distance $r$ from the first root vertex are mapped into vertical lines through $(r, 0)$, the only possible condensation point in question being $\infty$. We leave it to the reader to verify uniqueness. Below, we shall consider rooted trees (finite or infinite) as planar maps via this correspondence.

We are now ready to describe the mapping $\mathcal{Q}$, following closely [12], Section 3.4. Let $\omega$ be a tree in $\mathcal{C}$, considered as a planar map with condensation point for sequences as in Assumption ($\alpha$) equal to $c$, which we assume is in $\mathbb{R}^2$ (as also indicated in the figures of this section). By $F_0$ we denote the complement of $\omega$ in $\mathbb{R}^2$. A *corner* of $F_0$ is a sector between two consecutive edges around a vertex. A vertex of degree $k$ defines $k$ corners. The label of a corner is by definition the label of the corresponding vertex. A labeled tree $\omega \in \mathcal{C}$ has a finite number $C_k(\omega) \geq N_k(\omega)$ of corners with label $k$.

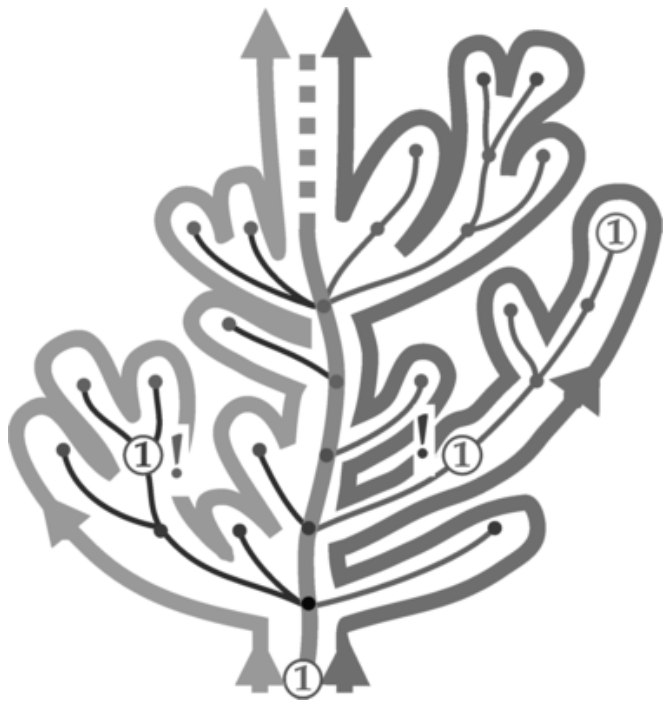

FIG. 3. *Contour traversal* (a) *of the left side and* (b) *of the right side of the spine. The last occurrence of label 1 is signaled by an exclamation mark.*

The image $\mathcal{Q}(\omega)$ is defined in three steps.

STEP 1. A vertex $v_0$ with label 0 is placed in $F_0 \setminus \{c\}$ and one edge is added between this vertex and each of the $C_1(\omega) < +\infty$ corners with label 1. Notice that this is possible because $\omega$ has only one spine. The new root is taken to be the edge that arrives from $v_0$ at the corner before the root of $\omega$.

After Step 1, a uniquely defined rooted planar map $\mathcal{M}_0$ with $C_1(\omega) - 1$ faces has been obtained [see Figure 4, with $C_1(\omega) = 7$]. It is natural to consider the complement of $\mathcal{M}_0$ and its faces as an additional face, which we shall call the *infinite face*: whereas $C_1(\omega) < +\infty$, there is a last occurrence of the label 1 when one does a contour traversal, as indicated on Figure 3, of the left (resp. right) side of the spine, and this corresponds to some corner $c_\ell$ (resp. $c_r$) with label 1. The infinite face is the one with corner $c_\ell$-$v_0$-$c_r$. The other faces are bounded by edges joining $v_0$ to two corners on the same side of the spine and are thus finite.

The next steps take place independently in each of those faces and will be described for a generic face $F$ of $\mathcal{M}_0$ (the infinite face does not require special treatment). Let $k$ be the degree of $F$ ($k$ can be infinite and by construction $k \geq 3$). Among the corners of $F$, only one belongs to $v_0$ and has label 0. If $F$ is finite, let the corners be numbered from 0 to $k-1$ in clockwise order along the border, starting with $v_0$. Otherwise, $F$ contains infinitely many corners both on the left and right side of the spine. Let the corners on the right of the spine be numbered by nonnegative integers, in clockwise order, starting with $v_0$, and let the corners on the left of the spine be numbered by negative integers, in counterclockwise order, starting right after $v_0$. Moreover, let $\ell(i)$ be the label of corner $i$ [so that $\ell(0) = 0$ and in the finite case $\ell(1) = \ell(k-1) = 1$, while in the infinite case $\ell(1) = \ell(-1) = 1$]. In Figure 5 the corners are explicitly represented together with their numbering for the infinite face.

STEP 2. In each face, the function successor $s$ is defined for all corners, except the corner at $v_0$, by

$$s(i) = \min\{j \triangleright i | \ell(j) = \ell(i) - 1\},$$

in which $j \triangleright i$ has the same meaning as $j > i$ for couples of positive integers and also for couples of negative integers, but we assume that a negative integer is larger than a positive integer. More precisely, if $j \leq 0$ and $i > 0$, then $j \triangleright i$.

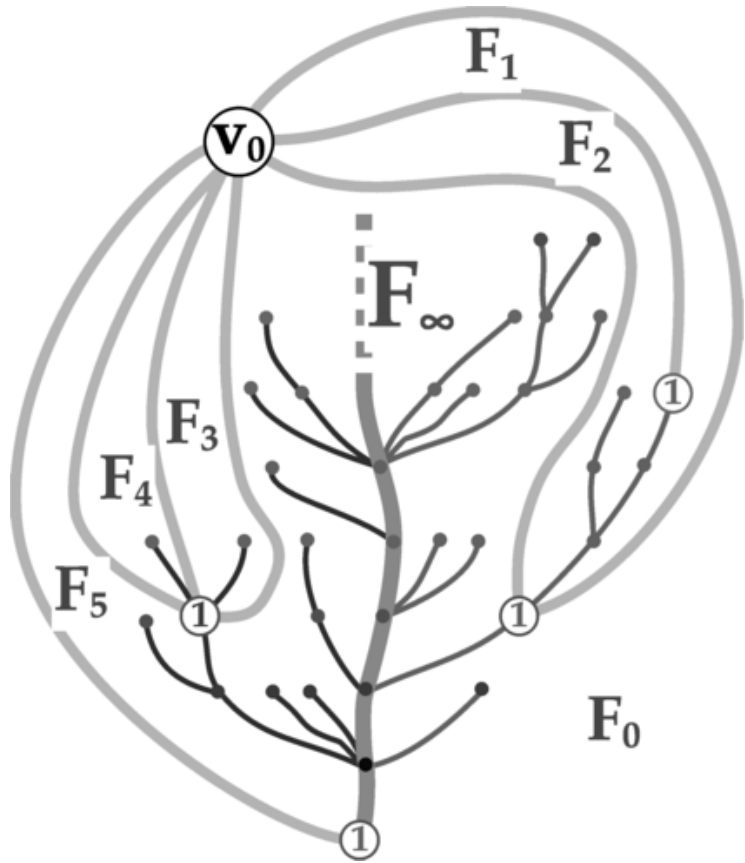

Fig. 4. *Step 1, leading to finite faces $F_i$, $0 \leq i \leq 5$, and the infinite face $F_\infty$.*

STEP 3. For each corner $i \geq 2$ such that $s(i) \neq i + 1$, a chord $(i, s(i))$ is added inside the face. This can be done in such a way that the various chords do not intersect (Property 6.1 below).

Once this construction has been carried out in each face, a planar map $\mathcal{M}'$ is obtained:

STEP 4. All edges of $\mathcal{M}'$ with the same label at both ends are deleted. The resulting map is a quadrangulation $\mathcal{M} = \mathcal{Q}(\omega)$ that satisfies Assumption $(\beta)$ (Property 6.2 below).

Note that a chord $(i, s(i))$ of the infinite face can very well join both sides of the spine $(is(i) < 0)$ if, in the contour traversal of the the right side of the spine, the corner $i$ appears after the last occurrence of label $\ell(i) - 1$ (see Figure 6). Let us first prove the two properties that validate the preceding construction.

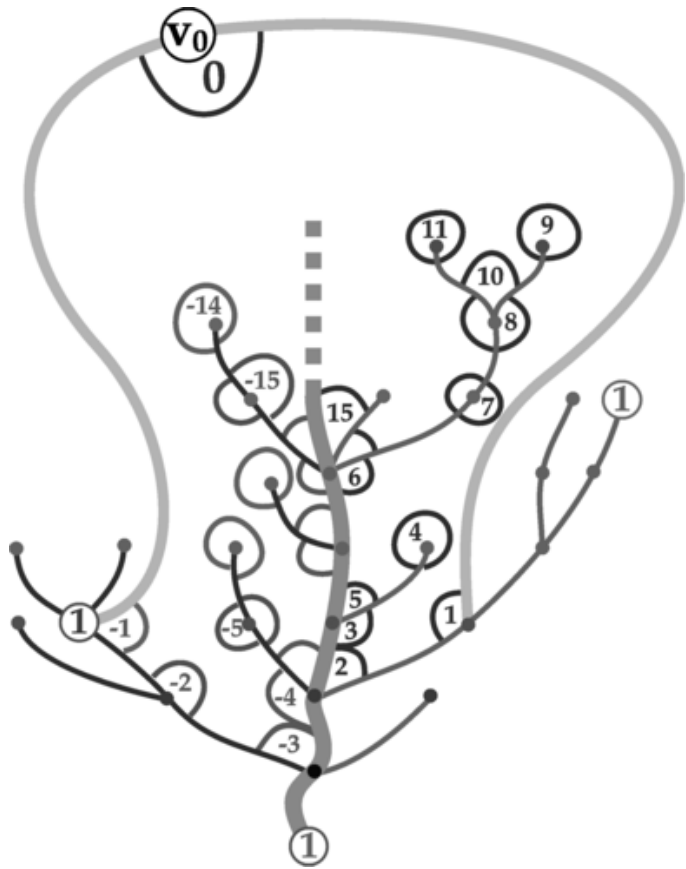

Fig. 5. *Numbering the corners of the infinite face.*

PROPERTY 6.1. *The chords $(i, s(i))$ do not intersect.*

PROOF. Suppose that two chords $(i, s(i))$ and $(j, s(j))$ cross each other. Perhaps upon exchanging $i$ and $j$, one has $i \triangleleft j \triangleleft s(i) \triangleleft s(j)$. The first two inequalities imply, together with the definition of $s$, that $\ell(j) > \ell(s(i))$, while the two last inequalities imply $\ell(s(i)) \geq \ell(j)$. This is a contradiction. □

PROPERTY 6.2. *The faces of $\mathcal{M}'$ are of one of the two types in Figure 7: either triangular with labels $e, e+1, e+1$ or quadrangular with labels $e, e+1, e+2, e+1$.*

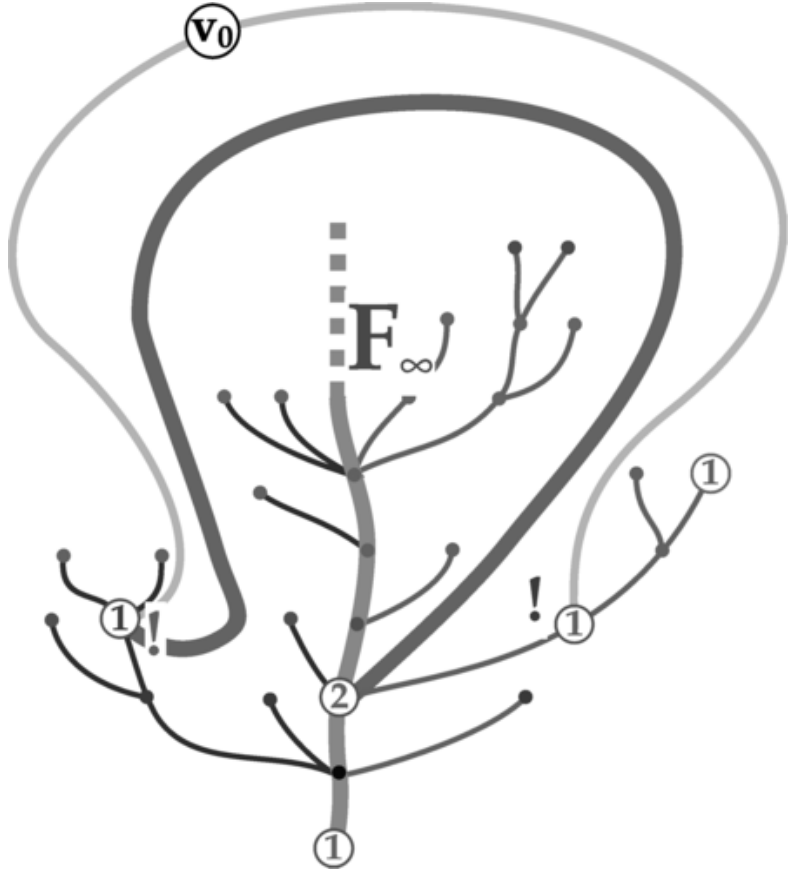

FIG. 6. *Inside $F_\infty$, a chord that joins two sides of the spine, with $\ell(2) = 2$, $\ell(-1) = 1$ and $s(2) = -1 \triangleright 2$.*

PROOF. Let $f$ be a face of $\mathcal{M}'$. Then $f$ is included in a face $F$ of $\mathcal{M}_0$ so that its corners inherit the numbering and labeling of those of $F$. Let $j$ be the corner with largest number (w.r.t. $\triangleleft$) in $f$, let $e = \ell(j)$ and let $i_1 \triangleleft i_2 \triangleleft j$ be the two neighbors of $j$ in $f$ (cf. Figure 7). The two latter corners both have label $e + 1$, because the edge $(i_1, j)$ has to be a chord so that $j = s(i_1)$ and, whereas $i_1 \triangleleft i_2 \triangleleft j$, this implies $\ell(i_2) \geq \ell(i_1)$ and hence $\ell(i_1) = \ell(i_2) = \ell(j) + 1$.

By construction, no other chord leaves $i_1$ so that the face is bordered by the edge $(i_1, i_1 + 1)$ of $F$. Either of two cases may occur:

- We have $i_1 + 1 = i_2$ and the face is triangular.
- The face is quadrangular: indeed the corner $i_1 + 1$ has label $e + 2$ [otherwise chord $(i_1 + 1, j)$ would exclude $i_2$ from the face] and the chord leaving $i_1 + 1$ goes to $i_2$ [otherwise $s(i_1 + 1) \neq i_2$, so that $\ell(s(i_1 + 1)) = \ell(i_1 + 1) - 1 = \ell(i_2)$ and chord $(s(i_1 + 1), j)$ would exclude $i_2$].

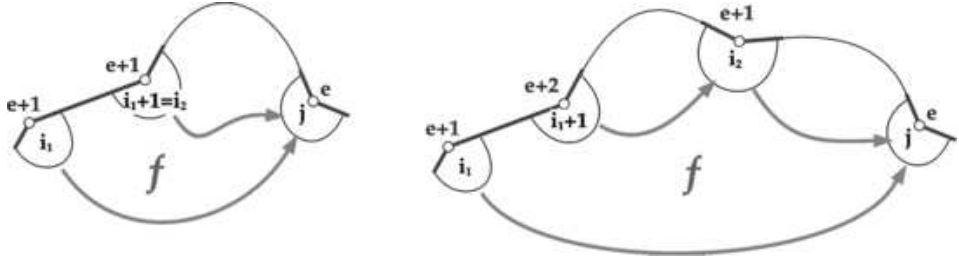

Fig. 7. *Two possible sizes for $f$: triangular or quadrangular.*

Observe finally that the deletion of edges with the same labels at both ends will join triangular faces pairwise to form quadrangular faces. □

PROPERTY 6.3. *The labels of the vertices in the tree $\omega$ are the distances to the root of the corresponding vertices in $\mathcal{M} = \mathcal{Q}(\omega)$.*

PROOF. By construction, the variation along the edges of the quadrangulation $\mathcal{Q}(\omega)$ is $0$ or $\pm 1$. Thus any path from a vertex $v$ with label $k$ to $v_0$ has length at least $k$. Choose a face $F$ of $\mathcal{M}_0$, such that $v$ is a corner of $F$ and let $s$ be the successor function associated with $F$. Then $v, s(v), s^2(v), \ldots, s^k(v)$ is a length-$k$ path in $\mathcal{Q}(\omega)$ that ends at $v_0$. □

6.3. *Properties of $\mathcal{Q}(\omega)$.* In this section, for the sake of brevity, $\mathcal{D}(\omega)$ will denote $\mathcal{D}(\mathcal{Q}(\omega))$. As remarked earlier, since by construction $\mathcal{Q}(\omega)$ satisfies Assumption ($\beta$), $\mathcal{D}(\omega)$ is an open connected subset of $S^2$. It should be noted that $\mathcal{D}(\omega)$ depends on the way the chords and edges in $\omega$ are drawn in the plane. However, it is easy to see that, as a planar map, $\mathcal{Q}(\omega)$ is unique by construction and that $\mathcal{Q}$ is an injective mapping. In particular, $\mathcal{D}(\omega)$ is unique up to homeomorphisms for any given $\omega \in \mathcal{C}$.

Let now $\bar{\mu}$ be the measure on $\mathcal{Q}(\mathcal{C})$ obtained by transporting $\mu$ from $\mathcal{C}$, that is,

$$\bar{\mu}(A) = \mu(\mathcal{Q}^{-1}(A))$$

for subsets $A \subseteq \mathcal{Q}(\mathcal{C})$, such that $\mathcal{Q}^{-1}(A)$ is $\mu$-measurable in $\mathcal{C}$. Also, let $B_r(\mathcal{M})$ denote the ball of radius $r$ in a quadrangulation $\mathcal{M}$, that is, the finite planar map whose vertices are those of $\mathcal{M}$ with (graph) distance from the first root vertex less than or equal to $r$, together with the edges that connect them (see Section 1). By Property 6.2 we can then reformulate Theorem 5.9 as follows.

THEOREM 6.4. *Let $|B_r(\mathcal{M})|$ denote the number of vertices in $B_r(\mathcal{M})$. Then*

$$\mathbb{E}_{\bar{\mu}}[|B_r|] = \Theta(r^4).$$

It should be noted that the measure $\bar{\mu}$ is, in fact, a weak limit of the uniform measures $\bar{\mu}_N, N \in \mathbb{N}$, on finite quadrangulations with $N$ quadrangles, for which reason we call $\bar{\mu}$ the *uniform probability measure on infinite quadrangulations.* To formulate this more precisely, note first that, as mentioned previously, the map $\mathcal{Q}$ extends to the subset

$$\mathcal{C}_0 = \mathcal{C} \cup \mathcal{W}$$

of $\overline{\mathcal{W}}$, obtained by adjoining the finite labeled trees to $\mathcal{C}$, and it maps $\mathcal{W}'_N$ bijectively onto the set of quadrangulations with $N$ quadrangles. Denoting the extended map also by $\mathcal{Q}$, we have that $\mathcal{Q}$ maps $\mathcal{C}_0$ bijectively onto its image $\mathcal{Q}(\mathcal{C}_0)$. Now define the metric $\bar{d}$ on $\mathcal{Q}(\mathcal{C}_0)$ such that $\mathcal{Q}$ is an isometry from $(\mathcal{C}_0, d)$ to $(\mathcal{Q}(\mathcal{C}_0), \bar{d})$. Then $\bar{\mu}_N, N \in \mathbb{N}$, and $\bar{\mu}$ can naturally be considered as Borel probability measures on $\mathcal{Q}(\mathcal{C}_0)$, which under $\mathcal{Q}$ correspond to the measures $\mu_N, N \in \mathbb{N}$, and $\mu$, respectively, regarded as measures on $\mathcal{C}_0$. Hence $\bar{\mu}_N \to \bar{\mu}$ on $\mathcal{Q}(\mathcal{C}_0)$ is equivalent to $\mu_N \to \mu$ on $\mathcal{C}_0$. However, the latter convergence follows from Theorem 2.2 in [9] by observing that the denumerable family

$$\mathcal{V}_0 = \{\mathcal{B}_{1/r}(\omega) \cap \mathcal{C}_0 | r \in \mathbb{N}, |\omega| < \infty\}$$

is closed under finite intersections and generates the open sets in $\mathcal{C}_0$ by the same arguments as in the proof of Theorem 3.1 and, moreover, we have proven there that $\mu_N(A) \to \mu(A)$ as $N \to \infty$ for all $A \in \mathcal{V}_0$.

We have thus established the following convergence result for the uniform measures on finite quadrangulations.

PROPOSITION 6.5. *The sequence $\bar{\mu}_N, N \in \mathbb{N}$, converges weakly to the uniform probability measure $\bar{\mu}$ on $\mathcal{Q}(\mathcal{C}_0)$.*

Our final results concern the shape of the domain $\mathcal{D}(\omega)$, which should be compared with Theorem 1.10 in [8].

THEOREM 6.6. *For any well-labeled tree $\omega \in \mathcal{C}$, the complement of $\mathcal{D}(\omega)$ in $S^2$ is connected. As such, the domain $\mathcal{D}(\omega)$ is homeomorphic to a disc and, in particular, it has exactly one boundary component.*

PROOF. According to [25], or [26], Chapter 13, only the first assertion of the theorem needs a proof: Assume on the contrary that the complement of $\mathcal{D}(\omega)$ is not connected, such that we can write

$$S^2 \setminus \mathcal{D}(\omega) = K_1 \cup K_2,$$

where $K_1$ and $K_2$ are nonempty compact subsets of $S^2$ contained in two disjoint open sets $O_1$ and $O_2$, respectively. From Assumption $(\alpha)$, it follows that there is only a finite number of edges in $\mathcal{Q}(\omega)$ that are not contained in $O_1 \cup O_2$, and clearly each of the sets $O_1$ and $O_2$ contains an infinite number of edges. Let now $R_0 > 0$ be large enough such that $B_{R_0}(\mathcal{Q}(\omega))$ contains all edges that are not contained in $O_1 \cup O_2$. It follows that if we remove from $\mathcal{D}(\omega)$ any finite set of faces that contain all vertices in $B_{R_0}(\omega)$, then the remaining part of $\mathcal{D}(\omega)$ is not connected. We shall now obtain a contradiction with this statement.

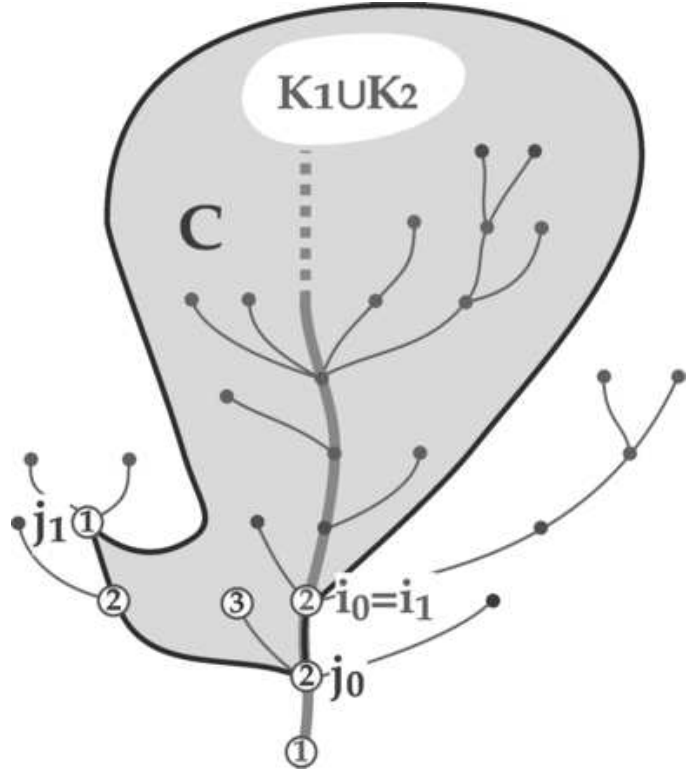

Fig. 8.   *The closed curve C for R = 2.*

Let $R \geq 2$ be a fixed integer. Let $i_1$ be the vertex with the last occurrence of label $R$ by (clockwise) contour traversal on the right hand side of $\omega$ and, similarly, let $j_1$ be the vertex with last occurrence of label $R-1$ by (counterclockwise) contour traversal on the left-hand side of $\omega$. Furthermore, let $i_0$ (resp. $j_0$) be the vertex on the spine of $\omega$ at which the branch that contains $i_1$ (resp. $j_1$) is attached. We then obtain a closed curve $C$ in $\mathcal{M}'$ made up of the shortest paths that connect $j_1$ to $j_0$, $j_0$ to $i_0$ and $i_0$ to $i_1$ together with the chord from $i_1$ to $j_1$ (see Figure 8). By construction, all vertices enclosed by this curve have labels larger than or equal to $R-1$ and only finitely many vertices are outside $C$. Obviously, the faces of $\mathcal{M}'$ enclosed by $C$ form a connected set. By adding to this the triangles that correspond to cutting edges in $C$ with equal labels, we obtain a connected set of faces in $\mathcal{M} = \mathcal{Q}(\omega)$ all of whose vertices have labels larger than or equal to $R-2$ and such that only finitely many vertices in $\mathcal{Q}(\omega)$ are not in this set. Since this holds for arbitrary $R \geq 2$, we have established the claimed contradiction, thus finishing the proof. $\square$

**Acknowledgments.** One of the authors (B. Durhuus) wishes to express his gratitude to Institut Élie Cartan for generous hospitality extended to him during his stay there in March and April 2003, where most of the present work was done. He also wishes to thank Thórdur Jónsson for helpful discussions at the early stages of this work. The authors thank Lucas Gerin and an anonymous referee for their valuable comments.

INSTITUT ÉLIE CARTAN
UNIVERSITÉ HENRI POINCARÉ–NANCY I
BP 239
54506 VANDŒUVRE-LÈS-NANCY
FRANCE
E-MAIL: chassain@iecn.u-nancy.fr
URL: www.iecn.u-nancy.fr/~chassain

DEPARTMENT OF MATHEMATICS AND MAPHYSTO
UNIVERSITY OF COPENHAGEN
UNIVERSITETSPARKEN 5
DK-2100 COPENHAGEN Ø
DANMARK
E-MAIL: durhuus@math.ku.dk
URL: www.math.ku.dk/~durhuus